%
%
%

\input amstex

\input amssym
\input amssym.def

\magnification 1200
\loadmsbm
\parindent 0 cm


\define\nl{\bigskip\item{}}
\define\snl{\smallskip\item{}}
\define\inspr #1{\parindent=20pt\bigskip\bf\item{#1}}
\define\iinspr #1{\parindent=27pt\bigskip\bf\item{#1}}
\define\ainspr #1{\parindent=24pt\bigskip\bf\item{#1}}
\define\aiinspr #1{\parindent=31pt\bigskip\bf\item{#1}}

\define\einspr{\parindent=0cm\bigskip}

\define\ot{\otimes}

\centerline{\bf The Larson-Sweedler theorem for weak multiplier Hopf algebras}

\bigskip\bigskip
\centerline{\it Byung-Jay Kahng \rm $^{(1)}$ and \it Alfons Van Daele \rm $^{(2)}$} 
\bigskip\bigskip\bigskip
{\bf Abstract} 
\nl 
The Larson-Sweedler theorem says that a finite-dimensional bialgebra with a faithful integral is a Hopf algebra [L-S]. The result has been generalized to finite-dimensional weak Hopf algebras by Vecserny\'es [Ve]. In this paper, we show that the result is still true for weak multiplier Hopf algebras. 
\snl
The notion of a weak multiplier bialgebra was introduced by B\"ohm, G\'omez-Torrecillas and L\'opez-Centella in [B-G-L]. In this note it is shown that a weak multiplier bialgebra with a regular and full coproduct is a regular weak multiplier Hopf algebra if there is a faithful set of integrals. Weak multiplier Hopf algebras are introduced and studied in [VD-W3]. Integrals on (regular) weak multiplier Hopf algebras are treated in [VD-W5].
\snl
This result is important for the development of the theory of locally compact quantum groupoids in the operator algebra setting, see [K-VD2] and [K-VD3]. Our treatment of this material is motivated by the prospect of such a theory.
\vskip 5 cm
Date: {\it 1 March 2016} (Version 2.0)
\vskip 1 cm
\hrule
\bigskip
\parindent 0.7 cm
\item{($1$)} Department of Mathematics and Statistics, Canisius College Buffalo, NY 14208, USA.  {\it E-mail}: kahngb\@canisius.edu
\item{($2$)} Department of Mathematics, University of Leuven, Celestijnenlaan 200B,
B-3001 Heverlee, Belgium. {\it E-mail}: Alfons.VanDaele\@wis.kuleuven.be

\parindent 0 cm

\newpage

\bf 0. Introduction \rm
\nl
In this paper we will prove the {\it Larson-Sweedler theorem} for {\it weak multiplier Hopf algebras}. 
\snl
In order to explain our result, in this introduction we first recall some basic notions and related results.
\snl
i) To begin with, there is the well-known notion of a {\it bialgebra} $(A,\Delta)$. In this case $A$ is an (associative) algebra with an identity and $\Delta$ is a coproduct on $A$. A coproduct (or comultiplication) is a unital homomorphism $\Delta:A\to A\ot A$ that satisfies coassociativity, i.e.\ so that 
$$(\Delta\ot\iota)\Delta(a)=(\iota\ot\Delta)\Delta(a)\tag"(0.1)"$$ 
for all $a\in A$, where we use $\iota$ for the identity map on $A$. It is also assumed that there is a counit. This is a homomorphism $\varepsilon$ from the algebra $A$ to the scalar field such that 
$$(\varepsilon\ot \iota)\Delta(a)=a
\qquad\quad\text{and}\qquad\quad
(\iota\ot\varepsilon)\Delta(a)=a \tag"(0.2)"$$
for all $a$. Because the unit in an algebra is unique if it exists and because the same is true for a counit, there is no need to include these objects in the notation for a bialgebra. For the notion of a bialgebra, there are several standard references. We will collect some of them later in the introduction in the separate item '{\it Various references}'.
\snl
An easy but important result says that the linear dual of a finite-dimensional bialgebra is again, in a natural way, a bialgebra. The result is no longer true if the algebra is not assumed to be finite-dimensional. Again see the standard references about bialgebras.
\snl
ii) Next we recall the notion of a {\it weak bialgebra}. In this case, it is assumed that $A$ is a unital algebra and that $\Delta$ is a coassociative homomorphism from $A$ to $A\ot A$ but it is {\it no longer} required that the coproduct $\Delta$ is a {\it unital} homomorphism. On the other hand, one still assumes that there is a counit $\varepsilon$ on $A$ satisfying (0.2), but it is {\it no longer} required to be a {\it homomorphism}. Instead some weaker conditions are assumed. We refer to the standard references about weak bialgebras, again see further in this introduction.
\snl
Also here, the dual of a finite-dimensional weak bialgebra is again a weak bialgebra.
\snl
iii) We can generalize the case of a bialgebra, as considered in item i), into another direction. Now it is {\it no longer} assumed that the algebra $A$ has {\it an identity}. This condition is replaced by the weaker requirement that the product is non-degenerate and that the algebra is idempotent. For the coproduct one now takes a homomorphism $\Delta$ from $A$ to the {\it multiplier algebra} $M(A\ot A)$. The homomorphism $\Delta$ is assumed to be non-degenerate and this implies that there is a unique extension to a unital homomorphism from $M(A)\to M(A\ot A)$. Also the homomorphisms $\Delta\ot \iota$ and $\iota\ot\Delta$ have unique extensions and therefore it makes sense to require again coassociativity as in (0.1). Still the counit is  a homomorphism $\varepsilon$ from $A$ to the field satisfying the conditions (0.2). In  order to make these conditions precise, also the maps $\varepsilon\ot\iota$ and $\iota\ot\varepsilon$ have to be extended to the multiplier algebra $M(A\ot A)$. This results in the notion of a {\it multiplier bialgebra}. This case was first considered in [J-V]. However more recently, the notion has been defined also in [B-G-L] (see Theorem 2.11 in that paper and a remark made after the proof). The approach is different from [J-V] and it is considered as a special case of a weak multiplier bialgebra (see further in item iv)). The basic notions (like that of a multiplier algebra) that are needed here are found in e.g.\ [VD1]. See also the other references further in this introduction.
\snl
There is no relevant duality result here in the finite-dimensional case, because if $(A,\Delta)$ is a multiplier bialgebra and if $A$ is finite-dimensional, then $A$ has a unit and we are back in the case of a bialgebra as  considered in item i).
\snl
iv) The last case is that of a {\it weak multiplier bialgebra} $(A,\Delta)$. The idea is clear. The axioms are weakened in both directions. So, roughly speaking, the algebra is not assumed to be unital and the coproduct is no longer assumed to be non-degenerate. This case is considered in [B-G-L]. Because this is the basic concept for the current paper, we will recall this notion more precisely later, in Section 1.
\nl
In all these four cases, one can require the {\it existence of an antipode}. In the first case, for a bialgebra, it is an anti-homomorphism $S$ from $A$ to itself satisfying the conditions
$$m(S\ot \iota)\Delta(a)=\varepsilon(a)1
\qquad\quad\text{and}\qquad\quad
m(\iota\ot S)\Delta(a)=\varepsilon(a)1 \tag"(0.3)"$$
where $m$ is the multiplication map, considered from $A\ot A$ to $A$. If a bialgebra has an antipode, it is called a {\it Hopf algebra}.
\snl
Also in the three other more general cases, one can assume the existence of an antipode, but  the notion has to be adapted for each of these cases. A weak bialgebra with an antipode is a weak Hopf algebra, a multiplier bialgebra with an antipode is a multiplier Hopf algebra and a weak multiplier bialgebra with an antipode is a weak multiplier Hopf algebra. References for these cases are again given further in this introduction.
\snl
In general, the antipode need not be a bijective map from $A$ to itself. If that is the case, we call these objects {\it regular}.
\snl
In the finite-dimensional case, if the bialgebra is a Hopf algebra, the dual bialgebra is again a Hopf algebra. Similarly for a finite-dimensional weak Hopf algebra. Again, this is no longer true in the infinite-dimensional case.
\nl
Next there is the notion of {\it integrals}. If $(A,\Delta)$ is a bialgebra as in item i), a non-zero linear functional $\varphi$ on $A$ is called a left integral and a non-zero linear functional $\psi$ on $A$ is called a right integral if
$$(\iota\ot\varphi)\Delta(a)=\varphi(a)1
\qquad\quad\text{and}\qquad\quad
(\psi\ot\iota)\Delta(a)=\psi(a)1$$
for all $a\in A$. Integrals are unique (up to a scalar) if they exist.
\snl
For a finite-dimensional Hopf algebra, we always have the existence of integrals. They are automatically faithful.
\snl
The notion of integrals also exists for the other cases but again, it has to be modified so as to fit with the other axioms. In the case of multiplier bialgebras, integrals are also unique. However, in the case of a weak bialgebra integrals in general are not unique and the same is true for weak multiplier bialgebras. 
\snl
Again, in the case of a finite-dimensional weak Hopf algebra, integrals always exist. A single integral is not necessarily faithful here, but there is a set of faithful integrals (cf.\ Definition 2.10 in Section 2 for a precise formulation of this property). A necessary and sufficient condition for the existence of a single faithful integral is that the underlying algebra has a faithful functional, i.e.\ that it is a Frobenius algebra. See Section 3 in [B-N-S]. 
\snl
With integrals, duality can be extended to the non finite-dimensional case in a certain sense. If $(A,\Delta)$ is a regular multiplier Hopf algebra with integrals, there exists a suitable subspace of the linear dual of $A$ that can be made in a natural way into a regular multiplier Hopf algebra with integrals. Biduality holds in the sense that the dual of the dual is canonically isomorphic with the original. A similar result holds for regular weak multiplier Hopf algebras. For references about this duality, see later under the item '{\it Various references}'.
\nl
Finally, we are ready to describe the content of the {\it Larson-Sweedler theorem} as it can be proven in all these cases. Roughly speaking, the result says that, if a bialgebra has (enough) integrals, then there is an antipode.
\snl
More precisely, for a bialgebra we have that it is a Hopf algebra if it has integrals. For a multiplier bialgebra we get a multiplier Hopf algebra if there is an integral. For weak bialgebras one needs enough integrals in order to get a weak Hopf  algebra and the same is true in the case of multiplier weak bialgebras. We will again give references for the first three cases later in the introduction. The last result is precisely the one that will be proven in this paper.
\nl
\it Content of the paper \rm
\nl
In {\it Section} 1 we  recall the notion of a weak multiplier bialgebra (as it has been introduced recently in [B-G-L]) and see how this generalizes multiplier bialgebras, weak bialgebras and of course also bialgebras. We will  formulate the basic results that are needed for the proof of the main theorem in this paper. We will not give proofs, we refer to [B-G-L].
\snl
In {\it Section} 2 we will use these results about weak multiplier bialgebras to prove the Larson-Sweedler theorem in this setting. So, we will show that for any weak multiplier bialgebra, if there are enough integrals, there is an antipode and so it is a weak multiplier Hopf algebra. On the other hand, we will make the arguments more or less independent from the results obtained in the first section. In fact, we will start with a weaker set of assumptions than what is obtained for weak multiplier bialgebras in the first section. We will e.g.\ not assume the existence of a counit here. The motivation for this lies in the hope to use this approach for the further development of a theory of locally compact quantum groupoids (where counits in general are not considered). See the discussion in Section 3.
\snl
Indeed, in {\it Section} 3 we will discuss the results obtained in Section 1 and Section 2 further and make some remarks about the importance of the results for further research, in particular for the link of the algebraic theory with the operator algebraic approach. As mentioned already, the treatment in Section 2 is such that it is already preparing for this generalization into the direction of the operator algebra setting.
\snl
In {\it Appendix} A we discuss the differences in notation and terminology in this paper with the existing literature. We provide some kind of {\it dictionary}. Recall that this paper is, in some sense, preparing for a generalization to the operator algebraic context and that the algebraists and the operator algebraists  often use different conventions and notations. 
\snl
In {\it Appendix} B we look at the notions of separability and the Frobenius property as they appear in this note. Here again we also  discuss the relation of our conventions with the commonly used ones in the algebraic literature. 
\nl
\it Conventions and notations \rm
\nl
In this note, algebras are associative algebras over the field of complex numbers. They are in general not required to be unital, but the product is always assumed to be {\it non-degenerate} (as a bilinear form). We will also only work with {\it idempotent} algebras $A$, i.e.\ where $A=A^2$. The two conditions are fulfilled if the algebras have local units, but we will not need to assume that. It will follow from the other requirements. We use $A^{\text{op}}$ for the algebra obtained from $A$ by reversing the product.
\snl
For any non-degenerate algebra $A$ we can define its {\it multiplier algebra} $M(A)$. This can be characterized as the largest algebra containing $A$ as an essential ideal. If $A$ has an identity, then $A=M(A)$. We will also consider $A\ot A$ and its multiplier algebra $M(A\ot A)$. We have natural embeddings
$$A\ot A\subseteq M(A)\ot M(A)\subseteq M(A\ot A)$$
and in general, these inclusions are strict. 
\snl
We will use $1$ for the identity in $M(A)$. On the other hand, we use $\iota$ for the identity map from $A$ to itself. 
\snl
A coproduct $\Delta$ on $A$ is a homomorphism from $A$ to the multiplier algebra $M(A\ot A)$ satisfying certain properties. We will use $\Delta^{\text{cop}}$ for the coproduct on $A$ that is obtained by composing $\Delta$ with the flip map (extended to the multiplier algebra). For a coproduct, we will sometimes use the Sweedler notation and write $\Delta(a)=\sum_{(a)}a_{(1)}\ot a_{(2)}$ etc. This has to be done with care and it will be explained when we use this. There is also a new note devoted to the use of the Sweedler notation for multiplier Hopf algebras and weak multiplier Hopf algebras. See reference [VD6].
\snl
We will use the leg-numbering notation. If e.g.\ $E$ is an element in $M(A\ot A)$ we will use $E_{12}$ and $E_{23}$ for the elements $E\ot 1$ and $1\ot E$ in $M(A\ot A\ot A)$. Also we will use $E_{13}$ for $E$ as 'sitting' in this multiplier $M(A\ot A\ot A)$ with its 'first leg' in the first factor and with its 'second leg' in the third factor. More precisely $E_{13}=\zeta_{12}E_{23}$ where $\zeta_{12}$ flips the first two factors in $A\ot A\ot A$ and is extended to the multiplier algebra. 
\snl
Given a coproduct $\Delta$ on an algebra $A$ we can define four canonical maps $T_1$, $T_2$, $T_3$ and $T_4$ from $A\ot A$ to $M(A\ot A)$. The first one e.g.\ is given by 
$$T_1(a\ot b)=\Delta(a)(1\ot b)$$
when $a,b\in A$. For the other ones, we refer to Section 2 (see Notations 2.8). It is one of the requirements of the coproduct as used in this paper that these maps have ranges in $A\ot A$. We will use $\text{Ker}(T_i)$ and $\text{Ran}(T_i)$ for the kernel and the range of these maps.
\snl
Finally, see also the material about conventions and notations as discussed in the {\it dictionary} that we have provided in Appendix A of the paper.
\nl
\it Various references \rm
\nl 
For the theory of bialgebras and Hopf algebras, we have the standard references [Abe] and [Sw]. For the theory of weak bialgebras and weak Hopf algebras, we have [B-N-S] and [B-S]. See also [Ni], [N-V1], [N-V2] and [Va]. For the theory of multiplier Hopf algebras there is [VD1] and for integrals on multiplier Hopf algebras, the basic reference is [VD2]. Some information about the use of the Sweedler notation for multiplier Hopf algebras can be found in [VD3] and in the more recent work [VD6]. Multiplier bialgebras have been considered in [J-V].
\snl
The theory of weak multiplier Hopf algebras is developed in a series of papers [VD-W2], [VD-W3], [VD-W4] and [VD-W5]. Remark that [VD-W4] has two versions, [VD-W4.v1] and [VD-W4.v2]. The second version is more general in the sense that more results are included about not necessarily regular weak multiplier Hopf algebras.
We refer to [VD4] for a study of separability idempotents in relation with the theory of weak multiplier Hopf algebras. Also here there are two versions, [VD4.v1] and [VD4.v2]. As for [VD-W4] also for this paper, the second version includes the more general non-regular case. For the use of the Sweedler notation in this context, we refer to a paper that is currently still under preparation but should soon be available on the arXiv (cf.\ [VD6]).
\snl
The Larson-Sweedler theorem for Hopf algebras is found in the original paper [L-S]. For multiplier Hopf algebras, this result is obtained in [VD-W1]. The theorem for weak Hopf algebras is proven in [Ve]. In a forthcoming paper [VD7], some more reflections on the Larson-Sweedler theorem are found.
\snl
Some general references about separable algebras and Frobenius algebras are [DM-I] and [Abr]. A categorical approach is found in [St]. In the theory of weak (multiplier) Hopf algebras, these two concepts are connected and possible references then are [Sz], [Sc] and [K-S].
\nl
\bf Acknowledgments \rm
\nl
The first-named author (Byung-Jay Kahng) is especially grateful to his coauthor (Alfons Van Daele), for being a helpful mentor to him.  The current work was 
mostly done during his sabbatical leave visit to the University of Leuven during 2012/2013.  He is very much grateful to both Alfons Van Daele and the mathematics department at University of Leuven for their warm hospitality. He also wishes to thank Michel Enock for inspiring discussions on the topological version of the theory, on measured quantum groupoids. Finally, He wishes to thank the U.S.\ Fulbright Foundation for providing the financial support for his visit to Belgium.  
\snl
The second-named author (Alfons Van Daele) thanks Gabriella B\"ohm for the hospitality while he was visiting Budapest where part of this work was done and for private communications further on this subject. He is also grateful to Gabriella B\"ohm, Jos\'e G\'omez-Torrecillas and Esperanza L\'opez-Centella for providing a preliminary version of their work on weak multiplier bialgebras and for fruitful discussions about this work.
\nl\nl

\bf 1. Weak multiplier bialgebras \rm
\nl
First in this section we recall the notion of a {\it weak multiplier bialgebra} as introduced by B\"ohm, G\'omez-Torrecillas and L\'opez-Centella in [B-G-L]. We will also compare this notion here with the other (more restrictive) cases  that we have formulated in the introduction. We will give the main results about these objects, needed for the proof of the Larson-Sweedler theorem as we will treat it in the next section. 
\snl
We will not include proofs. We refer to [B-G-L]. It should be mentioned however that our set of assumptions is slightly different from the ones used in [B-G-L]. Fortunately, this makes little difference for the proofs. We will make more comments on this difference later, in Section 2.
\nl
\it The concept of a weak multiplier bialgebra \rm
\nl
A {\it weak multiplier bialgebra} is a pair $(A,\Delta)$ of an algebra $A$ with a coproduct $\Delta$ satisfying certain properties.
\snl
The algebra need not be unital, but is assumed to be non-degenerate in the sense that the product as a bilinear map is non-degenerate. It is also assumed to be idempotent. This means that any element of $A$ is the sum of products of two elements in $A$. The condition is written as $A=A^2$. These two conditions are automatic if the algebra is unital. They are also automatically fulfilled if the algebra has local units. In fact, as we will see further, the other conditions imposed on the pair $(A,\Delta)$ will imply that $A$ has to  have local units.
\snl
Recall also from the introduction that we only work with algebras over the field of complex numbers.
\snl
By a {\it coproduct} on $A$ we mean a homomorphism  $\Delta$ from $A$ to the multiplier algebra $M(A\ot A)$. It is assumed that all four elements of the form
$$\align
&\Delta(a)(1\ot b) \qquad\qquad\qquad (c\ot 1)\Delta(a) \tag"(1.1)"\\
&\Delta(a)(b\ot 1) \qquad\qquad\qquad (1\ot c)\Delta(a) \tag"(1.2)"
\endalign$$
belong to $A\ot A$ for all $a,b,c\in A$. Moreover $\Delta$ has to be coassociative in the sense that e.g.
$$(c\ot 1\ot 1)(\Delta\ot\iota)(\Delta(a)(1\ot b))=(\iota\ot\Delta)((c\ot 1)\Delta(a))(1\ot 1\ot b) \tag"(1.3)"$$
is true for all $a,b,c\in A$.
\snl
Remark that we use (1.1) to formulate coassociativity as in (1.3). In general, for a coproduct on an algebra $A$ that is not assumed to be unital, only (1.1) and not (1.2) is required. If also (1.2) is assumed, the coproduct is called {\it regular} (see e.g.\ [VD-W3]). So, in this paper, we assume from the very beginning that the coproduct is regular. This is not done in [B-G-L] while on the other hand, in order to obtain the main results, eventually also in [B-G-L] the coproduct is assumed to be regular. We will give some comments in Section 3 (Conclusions and further research) on the possibility to drop regularity as a condition for the results in this paper.
\snl
For a regular coproduct, it is possible to formulate coassociativity with the factors on the other side as
$$(\Delta\ot\iota)((1\ot b)\Delta(a))(c\ot 1\ot 1)=(1\ot 1\ot b)(\iota\ot\Delta)(\Delta(a)(c\ot 1)) \tag"(1.4)"$$
for all $a,b,c\in A$. Non-degeneracy of the product in $A$ will imply that (1.4) for all $a,b,c$ is equivalent with (1.3) for all $a,b,c$. 
\snl
Of course, if the algebra has a unit, these extra conditions are automatic and we just have the original condition of coassociativity. For algebras without a unit, this is a little more subtle. See e.g.\ a recent note on coassociativity for coproducts on algebras without identity [VD8].
\snl
The coproduct is also assumed to be {\it full}. Roughly speaking, this means that the legs of $\Delta(A)$ are all of $A$. More precisely, if $V,W$ are subspaces of $A$ so that
$$\Delta(A)(1\ot A)\subseteq V\ot A
\qquad\quad\text{and}\qquad\quad
(A\ot 1)\Delta(A)\subseteq A\ot W,$$
then $V=A$ and $W=A$. Because we assume the product to be regular, we can also express these conditions by multiplying on the other side. It can be shown that for a full coproduct, the span of elements of the form $(\omega\ot \iota)((a\ot 1)\Delta(b))$ where $a,b\in A$ and $\omega$ is a linear functional on $A$ is all of $A$. The same is true when we put $a$ on the right hand side, and if we act on the second leg. For more information about this notion and this property, we refer to [VD-W1] and also e.g.\ [VD-W2].
\snl
Again, if we compare with the original work of [B-G-L], we assume fullness of the product from the very beginning in this note. In [B-G-L] it is eventually assumed to obtain the main results. See e.g.\ Theorem 3.13 in [B-G-L].
\snl
For the coproduct we also have the following result. The proof is straightforward and can be found e.g.\ in Lemma 3.3 of [VD-W2]. 

\inspr{1.1} Lemma \rm If there exists an idempotent element $E\in M(A\ot A)$ satisfying
$$\Delta(A)(A\ot A)=E(A\ot A) 
\qquad\quad\text{and}\qquad\quad
(A\ot A)\Delta(A)=(A\ot A)E, \tag"(1.5)"$$ 
this idempotent is unique. It is then the smallest idempotent $E$ in $M(A\ot A)$ satisfying 
$$E\Delta(a)=\Delta(a)
\qquad\quad\text{and}\qquad\quad
\Delta(a)E=\Delta(a)$$
for all $a\in A$.
\hfill$\square$
\einspr
If this idempotent $E$ exists, it is usually referred to as {\it the canonical idempotent}.  
\snl
If it exists, it can be shown that the coproduct $\Delta$ has a unique extension to a homomorphism $\widetilde\Delta:M(A)\to M(A\ot A)$ provided we require that still
$$\widetilde\Delta(m)=E\widetilde\Delta(m)=\widetilde\Delta(m)E$$ 
for all $m\in M(A)$. In what follows we will denote this extension also by $\Delta$. Similarly, also the homomorphisms $\Delta\ot\iota$ and $\iota\ot\Delta$ have unique extensions to $M(A\ot A)$ requiring now, using the same symbols for these extensions, that
$$(\Delta\ot\iota)(m)=(E\ot 1)((\Delta\ot\iota)(m))=((\Delta\ot\iota)(m))(E\ot 1)
$$
whenever now $m\in M(A\ot A)$. Similarly for $\iota\ot\Delta$. Observe that we automatically have $(\Delta\ot\iota)E=(\iota\ot\Delta)E$.
\snl
For more information about this extension procedure and for the proofs of these properties,  we refer to the appendix in [VD-W2]. See also the recent note on coassociativity [VD8]. 
\snl
Using this result, we can formulate the following properties.

\inspr{1.2} Definition \rm The coproduct $\Delta$ is called {\it weakly non-degenerate} if an idempotent $E\in M(A\ot A)$ exist, satisfying the conditions (1.5) of the lemma. We also consider the condition 
$$(\Delta\ot\iota)E=(E\ot1)(1\ot E)=(1\ot E)(E\ot1).\tag"(1.6)"$$
and we refer to it as {\it weak comultiplicativity of the unit}.  
\hfill$\square$
\einspr

The above terminology needs to be explained. This is done in the following remark.

\inspr{1.3} Remark \rm
First, the term {\it weak non-degeneracy} of the coproduct is motivated by the theory of multiplier Hopf algebras where the condition of Lemma 1.1, with $E=1$, is indeed non-degeneracy of the coproduct. 
\snl
In the case of a unit, where $E=\Delta(1)$, the property (1.6) is called the  weak comultiplicativity of the unit in [B-N-S] (see Definition 2.1 in that paper). And although there is no unit in our case, we will still refer to this condition as {\it weak comultiplicativity of the unit}. 
\hfill$\square$
\einspr

We will give more comments on the last notion later (see a remark just before Definition 1.5).
\nl
Finally we come to the existence and requirements about the {\it counit}.
\snl
A counit is a linear map from $A$ to $\Bbb C$ satisfying
$$(\varepsilon\ot\iota)(\Delta(a)(1\ot b))=ab
\qquad\quad\text{and}\qquad\quad
(\iota\ot \varepsilon)((a\ot 1)\Delta(b))=ab$$
for all $a,b\in B$.
Because we assume that the coproduct is regular, we can equivalently require that
$$(\varepsilon\ot\iota)((1\ot a)\Delta(b))=ab
\qquad\quad\text{and}\qquad\quad
(\iota\ot \varepsilon)(\Delta(a)(1\ot b))=ab$$
for all $a,b\in B$. Because the coproduct is assumed to be full, the counit is unique. See e.g.\ Proposition 1.12 in [VD-W2].
\snl
Now, it is not assumed that the counit is multiplicative but only weakly multiplicative as in the following definition. 

\inspr{1.4} Definition \rm
The counit is called {\it weakly multiplicative} if
$$\align
\varepsilon(abc)&=(\varepsilon\ot\varepsilon)((a\ot 1)\Delta(b)(1\ot c))\\
\varepsilon(abc)&=(\varepsilon\ot\varepsilon)((1\ot a)\Delta(b)(c\ot 1))
\endalign$$
for all $a,b$.
\hfill$\square$\einspr

Remark that the conditions for weak multiplicativity of the counit are dual to the conditions of weak comultiplicativity of the unit as formulated above in Definition 1.2. To make such a statement precise, one can e.g.\ consider the finite-dimensional case, two algebras with a non-degenerate pairing and coproducts induced by the products. Then it is not hard to verify that the conditions of Definition 1.4 are indeed precisely the dual forms of the condition (1.6).
\snl
Note that the weak multiplicativity of the counit as formulated above, is generalizing the same notion in the case of a weak bialgebra, see Definition 2.1 in [B-N-S]. And because of this duality, it is indeed also justified to call the property of the coproduct (1.6) as assumed in the Definition 1.2  {\it weak comultiplicativity of the unit}.
\snl
Finally remark that the second condition in Definition 1.4 above can only be formulated for a regular coproduct. We refer again to the discussion in the last section of this paper where we consider possible future research on this subject.
\nl
Now we are ready for the {\it main definition} in this section.

\inspr{1.5} Definition \rm [B-G-L]
A {\it weak multiplier bialgebra} is a pair $(A,\Delta)$ of a non-degenerate idempotent algebra $A$ with a full coproduct $\Delta$ that is weakly non-degenerate in the sense of Definition 1.2. Weak comultiplicativity of the unit is assumed. There also is a counit that is assumed to be weakly multiplicative in the sense of Definition 1.4.
\hfill$\square$\einspr

Before we formulate the main results, let us make some more comments about the concept of a weak multiplier bialgebra as just formulated in the above definition. Also observe that the definition as given above, is more restrictive than the one originally given in [B-G-L].

\inspr{1.6} Remarks \rm
i) Assuming weak comultiplicativity of the unit is a condition about the behavior of $\Delta$ on the legs of $E$, together with requiring that these legs commute. It is quite natural as is explained in Section 3 of [VD-W2].
\snl
ii) As mentioned already, the assumptions about the weak multiplicativity of the counit are dual to the assumptions about the behavior of the coproduct on the legs of $E$ and therefore  also these are quite natural. 
\snl
iii) Finally, the assumptions on the counit as formulated in Definition 1.4 are slightly more general then the ones originally assumed by  B\"ohm, G\'omez-Torrecillas and  L\'opez-Centella in [B-G-L] but the difference is not essential. Remark however that we assume regularity of the coproduct from the vey beginning.
\hfill$\square$\einspr

Before we continue with formulating some of the main properties of weak multiplier bialgebras as obtained in [B-G-L], let us compare with the notions of bialgebras, multiplier bialgebras and weak bialgebras. 
\snl
First consider a bialgebra. In this case, the algebra is assumed to have an identity, $\Delta$ is assumed to be unital and we have $E=1$. So, obviously, the conditions in Definition 1.2 are satisfied. In this case, the counit is assumed to be a homomorphism and then also weak multiplicativity as in Definition 1.4 is true. In the case of a weak bialgebra, still the algebra $A$ is unital and $E=\Delta(1)$. The conditions in Definition 1.2 and for the counit in Definition 1.4 are now part of the axioms for a weak bialgebra (see e.g.\ [B-N-S]). The comparison with the case of a multiplier bialgebra is clear if we consider the notion as in [B-G-L] but it is not so obvious with the approach in [J-V], where this is also considered, but in a somewhat different way from the other cases.
\nl
\it Properties of weak multiplier bialgebras \rm
\nl
In what follows, $(A,\Delta)$ is a weak multiplier bialgebra as in Definition 1.5, $E$ is the canonical idempotent in $M(A\ot A)$ and $\varepsilon$ is the counit.
\snl
We first introduce the associated counital maps. 

\inspr{1.7} Proposition \rm
There exist four linear maps from $A$ to $M(A)$, given by
$$\align 
\varepsilon_s(a)&=(\iota\ot\varepsilon)((1\ot a)E)
\qquad\quad\text{and}\qquad\quad
\varepsilon_s'(a)=(\iota\ot\varepsilon)(E(1\ot a))\\
\varepsilon_t(a)&=(\varepsilon\ot\iota)(E(a\ot 1))
\qquad\quad\text{and}\qquad\quad
\varepsilon_t'(a)=(\varepsilon\ot\iota)((a\ot 1)E).
\endalign$$
\vskip -0.7cm \hfill$\square$
\einspr

It is clear that e.g.\ $\varepsilon_s(a)$ is defined as a right multiplier of $A$. However, it can be shown that it is actually an element in $M(A)$. In fact, one has that $(1\ot a)E$ belongs to $M(A)\ot A$ and from this, it would follow that $\varepsilon_s(a)\in M(A)$. Similarly for the other cases. These results however are not obvious and it is somewhat remarkable that weak multiplicativity of the counit is used to prove this.
\snl
Our notation is different from what is found in the original literature on weak bialgebras and weak Hopf algebras. We use the conventions as in [VD-W3]. For convenience of the reader, we have included a short appendix (Appendix A) to this paper where we give a {\it dictionary}, relating various notations in our paper with the notations as they are used e.g.\ in [B-N-S] and in [B-G-L].

\inspr{1.8} Proposition \rm
The ranges of the maps $\varepsilon_s$ and $\varepsilon_s'$ coincide. Similarly, the ranges of $\varepsilon_t$ and $\varepsilon_t'$ coincide. Moreover $\varepsilon_s(A)$  and $\varepsilon_t(A)$  are non-degenerate subalgebras of $M(A)$ and  their multiplier algebras
$M(\varepsilon_s(A))$ and $M(\varepsilon_t(A))$ can be considered as sitting in $M(A)$.
\hfill$\square$
\einspr 

In fact, we have the following characterization of the multiplier algebras.

\inspr{1.9} Proposition \rm
Let $A_s=M(\varepsilon_s(A))$ and $A_t=M(\varepsilon_t(A))$. Then we have
$$\align 
A_s&=\{x\in M(A) \mid \Delta(x)=E(1\ot x)\} \\
A_t&=\{y\in M(A) \mid \Delta(y)=(y\ot 1)E\}.
\endalign$$
\vskip -0.8cm \hfill$\square$
\einspr

It can be shown that, in some sense, $\varepsilon_s(A)$ is the left leg of $E$ while $\varepsilon_t(A)$ is the right leg of $E$. In particular {\it these algebras commute} in $M(A)$ and the same is still true for the multiplier algebras $A_s$ and $A_t$. This implies e.g.\ that we have similar formulas as in the proposition with the opposite product.
\snl
The main result is the following. It is the one we will consider in the next section.
We use $B$ for the algebra $\varepsilon_s(A)$ and $C$ for the algebra $\varepsilon_t(A)$.

\iinspr{1.10} Theorem \rm [B-G-L] The canonical idempotent $E$ belongs to the multiplier algebra \newline $M(B\ot C)$. It is a separability idempotent in the sense of [VD4.v1]. The canonical anti-isomorphisms $S_B:B\to C$ and $S_C:C\to B$, characterized by the formulas
$$E(x\ot 1)=E(1\ot S_B(x))
\qquad\quad\text{and}\qquad\quad
(1\ot y)E=(S_C(y)\ot 1)E$$
for $x\in B$ and $y\in C$ are given by 
$$S_B(\varepsilon_s'(a))=\varepsilon_t(a)
\qquad\quad\text{and}\qquad\quad
S_C(\varepsilon_t'(a))=\varepsilon_s(a)$$
for all $a\in A$.
\hfill$\square$\einspr

Because we have assumed from the very beginning that we are in the regular case, the separability element is also regular. Remark that in [VD4.v1], only regular separability idempotents are considered while in [V4.v2] also the non-regular case is studied. For this paper, this is irrelevant and so we essentially only need reference [VD4.v1]. 
\snl
Remark that it follows that the algebras $B$ and $C$ have local units. See Proposition 1.9 in [VD4.v1]. 
\nl
There exist what we call {\it left and right integrals} in [VD4.v1] and {\it distinguished linear functionals } in [VD4.v2]. We will adopt the better terminology of [VD4.v2].
\snl
These distinguished linear functionals are the linear functional $\varphi_B$ on $B$ and the linear functional $\varphi_C$ on $C$, characterized by the formulas
$$(\varphi_B\ot \iota)E=1
\qquad\quad\text{and}\qquad\quad
(\iota\ot\varphi_C)E=1
$$
and given by 
$$\varphi_B(\varepsilon_s(a))=\varepsilon(a)
\qquad\quad\text{and}\qquad\quad
\varphi_C(\varepsilon_t(a))=\varepsilon(a)
$$
for all $a\in A$.
\snl
The modular automorphisms of these functionals are precisely the compositions $S_BS_C$ on $C$ for $\varphi_C$ and $(S_CS_B)^{-1}$ on $B$ for $\varphi_B$. See Section 2 of [VD4.v1]. We will be using these functionals for the construction of the counit in Section 2.
\snl
Also for these notions, we are using a different terminology than what is common and again we refer to Appendix A and Appendix B for the relation with the more commonly used ones in the algebraic literature.
\snl
For further details about all this, in particular for the proofs, we refer to [B-G-L].
\nl
To finish, we just make a small remark about the involutive case. If $A$ is a $^*$-algebra, it is natural to assume that $\Delta$ is a $^*$-homomorphism. Regularity of the coproduct is automatic in this case. It will follow from the uniqueness of $E$ that it is self-adjoint. And from the uniqueness of the counit, it also follows that $\varepsilon$ is self-adjoint. Then it will not be difficult to discover the behavior of the counital maps obtained in Proposition 1.7 with respect to the involution. It will also follow that the ranges and their multiplier algebras $A_s$ and $A_t$ are $^*$-subalgebras as expected. For the antipodal maps, we find e.g.\ $S_C(S_B(x)^*)^*=x$ for all $x\in B$, also as expected from the general theory developed in [VD4.v1].
\nl\nl

\bf 2. The Larson-Sweedler theorem \rm
\nl
In this section, we start with a {\it non-degenerate idempotent algebra} $A$ and a {\it regular} and {\it full coproduct} $\Delta$ on $A$. We assume the existence of {\it an idempotent}  $E\in M(A\ot A)$ satisfying
$$\Delta(A)(A\ot A)=E(A\ot A)
\quad\quad\text{and}\quad\quad
(A\ot A)\Delta(A)=(A\ot A)E.$$
We know that such an idempotent is unique and it is the smallest one that satisfies 
$$E\Delta(a)=\Delta(a)E=\Delta(a)$$ 
for all $a\in A$.
\snl
Then we can extend $\Delta\ot\iota$ and $\iota\ot\Delta$ and we will have $(\Delta\ot\iota)E=(\iota\ot\Delta)E$ (see e.g. the appendix in [VD-W2] or the more recent note [VD8]). It is natural to assume also 
$$(\Delta\ot\iota)E=(1\ot E)(E\ot 1)=(E\ot 1)(1\ot E).\tag"(2.1)"$$
This has been motivated in [VD-W2].
\snl
Finally, we will assume that $E$ is a {\it regular separability idempotent} in the following sense.
\nl
First we {\it assume the  existence of two algebras $B$ and $C$} sitting in $M(A)$ in a non-degenerate way. By this we mean that $BA=AB=A$ and $CA=AC=A$. It follows that the products in $B$ and $C$ are non-degenerate. The embeddings extend to embeddings of $M(B)$ and $M(C)$ in $M(A)$. We will consider these multiplier algebras as subalgebras of $M(A)$. Similarly $B\ot C$ sits in $M(A\ot A)$ non-degenerately and we can consider $M(B\ot C)$ as sitting in $M(A\ot A)$ as well. 
We will denote $A_s=M(B)$ and $A_t=M(C)$. 
\snl
We now {\it assume further that} $E\in M(B\ot C)$ and that it is a regular {\it separability idempotent} in $M(B\ot C)$ in the sense of [VD4.v1]. We have collected the main definitions and results about this notion in Appendix B.
\snl
We can prove the following.

\inspr{2.1} Proposition \rm The algebras $B$ and $C$, if they exist, are completely determined by $E$.
\snl\bf Proof\rm:
Because $E$ is assumed to be a separability idempotent in $M(B\ot C)$, it follows that $C$ is the smallest subspace $V$ of $C$ so that $E(b\ot 1)\in B\ot V$ for all $b\in B$. We claim that it will also be the smallest subspace of $V$ of $M(A)$ so that $E(a\ot 1)\in A\ot V$ for all $a\in A$. This will imply that $C$ is uniquely determined by $E$ and the conditions imposed on $E$. Similarly for $B$.
\snl
To prove the claim, first observe that 
$$E(ba\ot 1)\in BA\ot C\subseteq A\ot C$$
for all $a\in A$. And as $BA=A$, it follows that $E(a\ot 1)\in A\ot C$ for all $a\in A$. 
\snl
On the other hand, assume that $V$ is a subspace of $M(A)$ so that $E(a\ot 1)\in A\ot V$ for all $a\in A$. We need to argue that $C\subseteq V$. To do this, assume that $\omega$ is a linear functional on $M(A)$ that is zero on elements of $V$. Because also 
$$(b\ot 1)E(a\ot 1)\in A\ot V$$
for all $b$ in $B$ and $a\in A$, we will have 
$(\iota\ot\omega)((b\ot 1)E(a\ot 1))=0$
for all $b\in B$ and $a\in A$. Because $(b\ot 1)E\in B\ot M(A)$, we can define 
$x=(\iota\ot\omega)((b\ot 1)E)$ as an element of $B$ and we see that $
xa=0$ for all $a$. This implies $x=0$. So $(\iota\ot\omega)((b\ot 1)E)=0$ for all $b$. This implies that $\omega=0$  on $C$ because $C$ is the right leg of $E$. 
\hfill$\square$
\einspr

Because of this result, it makes sense to say that we {\it require $E$ to be a separability idempotent} as we will do in the sequel of this section.
\snl
Before we continue however, we want to make one more remark about the 'legs' of $E$.

\inspr{2.2} Remark \rm A separability idempotent $E$ is assumed to be full. This means that the left leg of $E$ is all of $B$ and that the right leg of $C$ is all of $C$. This is when the legs are considered of $E$ being a separability idempotent in $M(B\ot C)$. In this setting however, we have $E\in M(A\ot A)$ and in the proof of the proposition, we see that $B$ will also be the left leg of $E$ as sitting in $M(A\ot A)$ and similarly for $C$. In this case, it means e.g.\ that elements of the form $(\iota\ot\omega(\cdot\,a))E$, where now $a\in A$ and $\omega$ is a linear functional on $A$, will be in $B$ and will span all of $B$. 
\hfill $\square$
\einspr

The reader is advised to compare with similar results about the legs of a regular coproduct as e.g.\ discussed in Section 1 of [VD-W1]. In fact, not only the results are similar, also the arguments are essentially the same. See also Definition 1.10 and Lemma 1.11 in [VD-W2] for some more information. 
\snl
In the following proposition, we will argue that the aforementioned conditions on $(A,\Delta)$ are fulfilled if it is a weak multiplier bialgebra (as defined and reviewed in the previous section). 

\inspr{2.3} Proposition \rm Let $(A,\Delta)$ be a weak multiplier bialgebra as in Definition 1.5 of the previous section. Then the associated idempotent $E$ (cf.\ Definition 1.2) is a regular separability idempotent as above.

\snl\bf Proof\rm: Of course we take $B=\varepsilon_s(A)$ and $C=\varepsilon_t(A)$. We have seen in Proposition 1.8 that these algebras sit in $M(A)$ in a non-degenerate way. In Theorem 1.10 it is stated that actually $E$ is in $M(\varepsilon_s(A)\ot\varepsilon_t(A))$ and that it is a separability idempotent.
\hfill $\square$
\einspr

The result is important and it shows that we could have started with the assumptions as in Section 1. On the other hand we do not from the very beginning assume that we have a counit for reasons we have explained already in the introduction. Indeed, since we are motivated by a possible development of the theory in an operator algebra framework, where the use of the counit is always avoided, we think it is important to start in this section {\it without} the assumption of the existence of a counit. See the discussion in Section 3.
\snl
In fact, we will show that a counit exists when we start with an idempotent $E$ as above and if we require the existence of sufficiently many integrals as we will do further.
\nl
We will use $S_B$  and $S_C$ for the canonical anti-isomorphisms from $B$ to $C$ and from $C$ to $B$ respectively defined by 
$$E(x\ot 1)=E(1\ot S_B(x))
\qquad\quad\text{and}\qquad\quad
(1\ot y)E=(S_C(y)\ot 1)E$$
for $x\in B$ and $y\in C$. See [VD4.v1] or Appendix B about separability in this paper. Observe that we assume regularity so that actually the maps $S_B$ and $S_C$ have range in $C$ and $B$, rather than in their multiplier algebras as in the non-regular case (as studied in [VD4.v2]).
\snl
For the multiplier algebras $A_s$ of $B$ and $A_t$ of $C$ we find the following property (obtained already in Proposition 1.9 in the case of a weak multiplier bialgebra).

\inspr{2.4} Proposition \rm
If $x\in A_s$ then $\Delta(x)=E(1\ot x)=(1\ot x)E$. If $y\in A_t$ then $\Delta(y)=(y\ot 1)E=E(y\ot 1)$. 
\snl\bf Proof\rm:
Using formula (2.1) and because of Remark 2.2, saying that the left leg of $E$ is $B$, it follows that for elements $x\in B$ we have $\Delta(x)=E(1\ot x)=(1\ot x)E$. Then this also follows for the elements in the multiplier algebra $A_s$ of $B$. Indeed, let $x\in M(B)$. If $b\in B$ we find
$$\Delta(b)\Delta(x)=\Delta(bx)=E(1\ot bx)=\Delta(b)(1\ot x).$$
Now we can multiply from the left with $\Delta(a)$ for any $a\in A$ and because it is assumed that $AB=A$, it follows that also
$\Delta(a)\Delta(x)=\Delta(a)(1\ot x)$
for all $a\in A$. We can now multiply with $p\ot q$ from the left for $p,q\in A$ and use that $(A\ot A)\Delta(A)=(A\ot A)E$. Then it follows that also $\Delta(x)=E(1\ot x)$. 
\snl
Similarly on the other side and for $C$ and its multiplier algebra $A_t$.
\hfill$\square$
\einspr

With these notions we can define {\it left} and {\it right integrals} on the pair $(A,\Delta)$.

\inspr{2.5} Definition \rm
A linear functional $\varphi$ on $A$ is called {\it left invariant} if $(\iota\ot\varphi)\Delta(a)\in A_t$ for all $a\in A$. A non-zero left invariant linear functional is called a {\it left integral}. Similarly a linear functional $\psi$ is called {\it right invariant} if $(\psi\ot\iota)\Delta(a)\in A_s$ for all $a$. If it is non-zero, it is called a {\it right integral}.
\hfill$\square$
\einspr

Because we assume that the coproduct is regular, the elements $(\iota\ot\varphi)\Delta(a)$ and $(\psi\ot\iota)\Delta(a)$ are well-defined in $M(A)$ for all $a\in A$ and so the conditions make sense.
\snl
Left and right integrals on regular weak multiplier Hopf algebras are studied in [VD-W5] (under preparation). We can not use the results from that paper anyway because we do not have all the properties of a weak multiplier Hopf algebra. So we will need to give some other arguments for the proofs of the following properties of the integrals. 
\nl
To begin with, we will need the following properties about the canonical idempotent $E$. They are found in Section 4 of [VD-W3] but as mentioned, we need to show that they are still valid under the weaker conditions in this setting. Remark that we are using the notations as in Proposition 4.7 in [VD-W3]. 
\snl
For the proof we refer to Appendix B because it is really a result about separability idempotents in general.

\inspr{2.6} Proposition \rm Denote 
$$\align F_1&=(\iota\ot S_C)E  \qquad\quad\text{and}\qquad\quad F_3=(\iota\ot S_B^{-1})E\tag"(2.2)"\\
	F_2&=(S_B\ot \iota)E \,\qquad\quad\text{and}\qquad\quad F_4=(S_C^{-1}\ot \iota)E.\tag"(2.3)"
\endalign$$
Then the elements $F_1$ and $F_3$ belong to $M(B\ot B)$ while $F_2$ and $F_4$ belong to $M(C\ot C)$. They satisfy 
$$\align E_{13}(F_1\ot 1) &= E_{13}(1\ot E) 
	\quad\quad\text{and}\quad\quad 
		(F_3\ot 1)E_{13}=(1\ot E)E_{13}\tag"(2.4)"\\
	(1\ot F_2)E_{13} &=(E\ot 1)E_{13}
	\quad\quad\text{and}\quad\quad 
	 	E_{13}(1\ot F_4)=E_{13}(E\ot 1).\tag"(2.5)"
\endalign$$
\vskip - 0.7cm \hfill$\square$
\einspr

Recall that the leg numbering notation, as used in this proposition, has been explained in the introduction.
\snl
We can now use these multipliers and obtain the following properties of the integrals.

\inspr{2.7} Proposition \rm If $\varphi$ is a left integral and $\psi$ a right integral we have for all $a$ in $A$ that
$$\align
(\iota\ot\varphi)\Delta(a)&=(\iota\ot\varphi)(F_2(1\ot a))\\
(\iota\ot\varphi)\Delta(a)&=(\iota\ot\varphi)((1\ot a)F_4)
\endalign$$
and
$$\align
(\psi\ot\iota)\Delta(a)&=(\psi\ot\iota)((a\ot 1)F_1)\\
(\psi\ot\iota)\Delta(a)&=(\psi\ot\iota)(F_3(a\ot 1))
\endalign$$

\snl\bf Proof\rm:
Consider a left integral $\varphi$. Take an element $a\in A$ and let $y=(\iota\ot\varphi)\Delta(a)$. By definition it is in $A_t$ and so 
$\Delta(y)=E(y\ot 1)$. This means that  
$$\Delta(y)=(\iota\ot\iota\ot\varphi)((E\ot 1)\Delta_{13}(a)).$$
Now we use that
$$(E\ot 1)E_{13}=(1\ot F_2)E_{13}$$
where $F_2$ is as in the previous proposition. This implies that also
$$(E\ot 1)\Delta_{13}(a)=(1\ot F_2)\Delta_{13}(a).$$
Then
$$\sum_{(a)} a_{(1)}\ot a_{(2)}\varphi(a_{(3)})=\Delta(y)
=(\iota\ot\iota\ot\varphi)((1\ot F_2)\Delta_{13}(a)).$$
This is true for all $a$ and we can use the fullness of the coproduct to obtain that then also
$$(\iota\ot\varphi)\Delta(a)=(\iota\ot\varphi)(F_2(1\ot a))$$
for all $a$. This proves the first formula.
\snl
If we start with the equality $\Delta(y)=(y\ot 1)E$ we arrive at the second formula in a completely similar way.
\snl
On the other hand, if $\psi$ is right invariant, and if we use that $x$, defined as $(\psi\ot\iota)\Delta(a)$ belongs to $A_s$, similar arguments as above will provide the two last formulas.
\hfill$\square$
\einspr

These formulas are found in [VD-W5]. They are also discussed in [VD7].
\snl
There is a nice consequence of these formulas. If follows that for any left integral $\varphi$ we have $\varphi(ya)=\varphi(a\sigma_C(y))$ for all $a\in A$ and $y\in A_t$ where $\sigma_C(y)=S_BS_C(y)$. Remark that $S_BS_C$ is the modular automorphism for the distinguished linear functional $\varphi_C$ on $E$ (see Proposition 2.3 in [VD4.v1] and the Appendix B). Indeed, we have $(S_CS_B\ot S_BS_C)E=E$ and therefore also $(\iota\ot S_BS_C)F_2=F_4$. 
Similarly we have for any right integral $\psi$ that $\psi(xa)=\psi(a\sigma_B(x))$ for all $a\in A$ and $x\in A_s$ where $\sigma_B(x)=S_B^{-1}S_C^{-1}(y)$. Again this is the modular automorphism for the unique right integral $\varphi_B$ on $E$ (see also Proposition 2.3 in [VD4.v1] and again also in the appendix). We will discuss this property again in Section 3 where we draw some conclusions and formulate further remarks.
\nl
Next recall the formulas for the {\it canonical maps} $T_1$, $T_2$, $T_3$ and $T_4$, associated with a regular coproduct (cf. Section 1 in [VD-W3]).

\inspr{2.8} Notations \rm
We denote
$$
\align T_1(a\ot b)&=\Delta(a)(1\ot b)
\qquad \quad\text{and}\qquad\quad
T_2(c\ot a)=(c\ot 1)\Delta(a)\\
T_3(a\ot b)&=(1\ot b)\Delta(a)
\qquad \quad\text{and}\qquad\quad
T_4(c\ot a)=\Delta(a)(c\ot 1)
\endalign$$
for all $a,b,c\in A$. Recall that the images of these maps are in $A\ot A$  because the coproduct is assumed to be regular.
\hfill$\square$
\einspr

Then we have the following result.

\inspr{2.9} Proposition \rm
Let $\varphi$ be a left integral and let $\psi$ be a right integral. Let $a,b\in A$.
We have
$$\align
&T_1((\iota\ot\iota\ot\varphi)(\Delta_{13}(a)\Delta_{23}(b)))=E(p\ot 1) 
\quad\qquad\text{where } 
p=(\iota\ot\varphi)(\Delta(a)(1\ot b))\\
&T_3((\iota\ot\iota\ot\varphi)(\Delta_{23}(a)\Delta_{13}(b)))=(p\ot 1)E 
\quad\qquad\text{where } 
p=(\iota\ot\varphi)((1\ot a)\Delta(b))
\endalign$$
and
$$\align
&T_2((\psi\ot\iota\ot\iota)(\Delta_{12}(a)\Delta_{13}(b)))=(1\ot q)E 
\quad\qquad\text{where } 
q=(\psi\ot\iota)((a\ot 1)\Delta(b))\\
&T_4((\psi\ot\iota\ot\iota)(\Delta_{13}(a)\Delta_{12}(b)))=E(1\ot q) 
\quad\qquad\text{where } 
q=(\psi\ot\iota)(\Delta(a)(b\ot 1)).
\endalign$$
These formulas have to be covered, i.e.\ we need to multiply with an element of the algebra, left or right to make these formulas meaningful. This is done in the proof below.

\snl\bf Proof\rm:
Take $a,b,c\in A$. First observe that 
$\Delta_{13}(a)\Delta_{23}(b)(1\ot c\ot 1)$
is well-defined as an element in the three fold tensor product $A\ot A\ot A$ by the regularity of the coproduct. Then we can apply the left integral $\varphi$ on the last factor to get elements in $A\ot A$. On such elements we can apply the canonical map $T_1$. 
Then we find
$$\align
T_1((\iota\ot\iota\ot\varphi)(\Delta_{13}(a)\Delta_{23}(b)&(1\ot c\ot 1)))\\
&\overset{\text{(a)}}\to{=}(\iota\ot\iota\ot\varphi)(((\iota\ot\Delta)\Delta(a))\Delta_{23}(b)(1\ot c\ot 1))\\
&=(\iota\ot\iota\ot\varphi)((\iota\ot\Delta)(\Delta(a)(1\ot b))(1\ot c\ot 1))\\
&\overset{\text{(b)}}\to{=}(\iota\ot\iota\ot\varphi)((1\ot F_2)(\Delta_{13}(a)(1\ot c\ot b)))\\
&\overset{\text{(c)}}\to{=}(\iota\ot\iota\ot\varphi)((E\ot 1)(\Delta_{13}(a)(1\ot c\ot b)))\\
&=E((\iota\ot\varphi)(\Delta(a)(1\ot b))\ot c)\\
&=E(p\ot c)
\endalign$$
with $p=(\iota\ot\varphi)(\Delta(a)(1\ot b))$. For the equality (a) we use co associativity of the coproduct, for (b) we use the results from Proposition 2.7 and for (c) we use Proposition 2.6.
\snl
The three other formulas are proven in the same way, using the other expressions in Proposition 2.7 and the formulas in Proposition 2.6.
\hfill$\square$
\einspr

Remark that we have used the formula $(E\ot 1)E_{13}=(1\ot F_2)E_{13}$ first to obtain the first formula in Proposition 2.7 and then again here in the proof above. The formulas in Proposition 2.9 do not explicitly involve these elements $F_i$. This seems to suggest that another proof is possible, not using these multipliers. We discuss this further in [VD7].
\nl
Now we assume the existence of {\it sufficiently many integrals} as follows.

\iinspr{2.10} Definition \rm 
We say that there is {\it a faithful set of left integrals} if, given $a\in A$, we must have $a=0$ if either $\varphi(ab)=0$ for all $b\in A$ and all left integrals or if $\varphi(ba)=0$ for all $b\in A$ and all left integrals on $A$. Similarly we define the notion of a {\it faithful set of right integrals}.
\hfill$\square$
\einspr

As a consequence of Proposition 2.9, we now get the following.

\iinspr{2.11} Proposition \rm
If there is a faithful set of {\it left} integrals, we find that 
$$T_1(A\ot A)=E(A\ot A)
\qquad\quad\text{and}\qquad\quad
T_3(A\ot A)=(A\ot A)E.$$
If there is a faithful set of {\it right} integrals, we find
$$T_2(A\ot A)=(A\ot A)E
\qquad\quad\text{and}\qquad\quad
T_4(A\ot A)=E(A\ot A).$$

\snl\bf Proof\rm:
Assume that there is a faithful set of left integrals. We claim that then $A$ is spanned by elements of the form
$$\{(\iota\ot\varphi)(\Delta(a)(1\ot b))\mid \varphi \text{ is a left integral, } a,b\in A\}.$$
Indeed, assume that there exists a linear functional $\omega$ on $A$ so that $\omega$ is $0$ on all such elements.
This will imply that for any left integral $\varphi$ and all elements $a,b,c\in A$ we have
$$\varphi((\omega\ot\iota)(\Delta(a)(1\ot b))c)=0.$$
As this is true for all $\varphi$ and all $c$, this implies that $(\omega\ot\iota)(\Delta(a)(1\ot b))=0$ for all $a,b$. Then the fullness of $\Delta$ will imply that $\omega=0$.
\snl
It is now an immediate consequence of the first formula in Proposition 2.9 that $T_1(A\ot A)=E(A\ot A)$.
\snl
A similar argument works for the three other cases.
\hfill$\square$
\einspr

This gives the necessary results for the {\it ranges} of the canonical maps. We now consider the {\it kernels} of these maps.

\iinspr{2.12} Proposition \rm
If there is a faithful set of {\it right} integrals, the kernels of the maps $T_1$ and $T_3$ are given by
$$
\text{Ker}(T_1)=(A\ot 1)(1-F_1)(1\ot A)
\quad\quad\text{and}\quad\quad
\text{Ker}(T_3)=(1\ot A)(1-F_3)(A\ot 1).$$
Similarly, if there is a faithful set of {\it left} integrals, the kernels of $T_2$ and $T_4$ are given by
$$
\text{Ker}(T_2)=(A\ot 1)(1-F_2)(1\ot A)
\quad\quad\text{and}\quad\quad
\text{Ker}(T_4)=(1\ot A)(1-F_4)(A\ot 1).$$

\snl\bf Proof\rm:
i) First we show that $T_1((a\ot 1)F_1(1\ot b))=T_1(a\ot b)$ for all $a,b$. Because $F_1=(\iota\ot S_C)E$ we have
$$(\Delta\ot\iota)(F_1)=(E\ot 1)(1\ot F_1).$$ 
Then we find
$$T_1((a\ot 1)F_1(1\ot b))=\Delta(a)E(1\ot m(F_1(1\ot b)))$$
where $m$ is multiplication form $A\ot A$ to $A$. Recall that $F_1=(\iota\ot S_C)E$ and that it is one of the properties of a separability idempotent that $mF_1=m(\iota\ot S_C)E=1$ (see [VD4.v1]). This implies that
$$T_1((a\ot 1)F_1(1\ot b))=\Delta(a)(1\ot b)$$
and this proves the claim. Remark that we do not need the integrals for this result.
\snl
ii) Conversely, assume that $\sum_i p_i\ot q_i\in \text{Ker}(T_1)$. This means that $\sum_i\Delta(p_i)(1\ot q_i)=0$. Multiply with $\Delta(a)$ from the left and apply a right invariant functional $\psi$ on the first leg. We find
$$\sum_i(\psi\ot\iota)(\Delta(ap_i)(1\ot q_i))=0$$
and by the formula in Proposition 2.7 we find that 
$$\sum_i(\psi\ot\iota)((ap_i\ot 1)F_1(1\ot q_i))=0.$$
As this holds for all $\psi$ and all $a$, it follows from the assumption that 
$$\sum_i(p_i\ot 1)F_1(1\ot q_i)=0.$$
This gives the right expression for the kernel of $T_1$. 
\snl
The other formulas are obtained in a similar way.
\hfill$\square$
\einspr

Remark that we need sufficiently many {\it left} integrals to find the right form for the range of $T_1$ and the kernel of $T_2$ whereas we need sufficiently many {\it right} integrals to get the right form for the kernel of $T_1$ and the range of $T_2$. Similarly for the maps $T_3$ and $T_4$.
\snl
We are almost ready with showing that $(A,\Delta)$ is a regular weak multiplier Hopf algebra. We just need to argue that the counit can also be constructed. This is done in the next proposition.

\iinspr{2.13} Proposition \rm
If $(A,\Delta)$ is as before and if there are enough left integrals, then there exists a counit. Similarly if there are enough right integrals.

\snl\bf Proof\rm: If $\varepsilon$ is a counit, we must have that 
$$\varepsilon((\iota\ot\varphi)(\Delta(a)(1\ot b))=\varphi(ab)\tag"(2.6)"$$
for all $a,b$ in $A$ and any linear functional, in particular any left integral $\varphi$ on $A$. We now try to define $\varepsilon$ with this formula. Therefore, let us now assume that we have elements $a_i,b_i$ in $A$ and left integrals $\varphi_i$ so that 
$$\sum_i(\iota\ot\varphi_i)(\Delta(a_i)(1\ot b_i))=0.\tag"(2.7)"$$ 
We need to show that $\sum_i\varphi_i(a_ib_i)=0$.
By Proposition 2.9 we find that 
$$\sum_i(\iota\ot\iota\ot\varphi_i)(\Delta_{13}(a_i)\Delta_{23}(b_i)(1\ot c\ot 1))$$
belongs to the kernel of $T_1$ for all $c$. It follows from Proposition 2.12 that 
$$\sum_i(\iota\ot\iota\ot\varphi_i)(\Delta_{13}(a_i)(F_1\ot 1)\Delta_{23}(b_i)(1\ot c\ot 1))=0.$$
We will now apply multiplication $m$ on this formula.  Using again that $mF_1=1$ and if we cancel $c$, we find
$\sum_i(\iota\ot\varphi_i)(\Delta(a_ib_i))=0$.
Using the first formula in Proposition 2.7 we find that 
$\sum_i(\iota\ot\varphi_i)(F_2(1\ot a_ib_i))=0$
and as $F_2=(S_B\ot \iota)E$ we get also
$\sum_i(\iota\ot\varphi_i)(E(1\ot a_ib_i))=0$.
And if now, we apply the distinguished linear functional $\varphi_C$ on $C$ satisfying 
$(\varphi_C\ot\iota)E=1$ we finally get 
$\sum_i\varphi_i(a_ib_i)=0$.
\snl
This proves that $\varepsilon$ is well-defined and that it satisfies (2.6). We need enough left integrals in order to have that all elements in $A$ are of the form as in (2.7). 
\snl
Then we have for all $a,b,b'$ and all $\varphi$ that
$$\varepsilon((\iota\ot\varphi)(\Delta(a)(1\ot bb')))=\varphi(abb')$$
and as this is true for all elements $b'$ and all left integrals $\varphi$ it follows that 
$$(\varepsilon\ot\iota)(\Delta(a)(1\ot b))=ab  \tag"(2.8)"$$
for all $a,b$. 
\snl
Now we can use coassociativity of $\Delta$ and the fullness of the coproduct to argue that also
$$(\iota\ot\varepsilon)((a\ot 1)\Delta(b))=ab$$
for all $a,b$. Indeed, take $a,b,c\in A$ and apply $\iota\ot\varepsilon\ot\iota$ on 
$$(\iota\ot\Delta)((a\ot 1)\Delta(b))(1\ot 1\ot c).$$
By (2.8) we find that this is equal to 
$(a\ot 1)\Delta(b)(1\ot c)$.
We now apply a linear functional $\omega$ on the last leg and use coassociativity. Because the coproduct is full, we can replace $(\iota\ot\omega)(\Delta(b)(1\ot c))$ by $b$ and  we find that 
$$(\iota\ot\varepsilon)((a\ot 1)\Delta(b))=ab.\tag"(2.9)"$$
\snl
Equation (2.8) and (2.9) show that $\varepsilon$ is the counit.
\snl
A similar proof works in the case of right integrals.
\hfill$\square$
\einspr

Combining all this, we can prove the following form of the Larson-Sweedler theorem for  weak multiplier Hopf algebras.

\iinspr{2.14} Theorem \rm
Assume that $(A,\Delta)$ is a pair of a non-degenerate idempotent algebra $A$ with a full and regular coproduct $\Delta$. Assume that there is an idempotent element $E$ in $M(A\ot A)$ satisfying the assumptions as formulated in the beginning of this section, in particular assumption (2.1). If there are a faithful set of left integrals and a faithful set of right integrals, then $(A,\Delta)$ is a regular weak multiplier Hopf algebra.

\snl\bf Proof\rm:
First we have to verify that the pair $(A,\Delta)$ satisfies the conditions of Definition 1.14 of [VD-W3]. 
\snl
We have a pair $(A,\Delta)$ of a non-degenerate idempotent algebra $A$ with a full coproduct $\Delta$ by assumption. We also have a counit as shown in Proposition 2.13.
\snl
There exist an idempotent $E\in M(A\ot A)$ giving the ranges of the canonical maps $T_1$ and $T_2$ as required. This is the content of Proposition 2.11. The idempotent is assumed to satisfy (2.1). This takes care of all the requirements in i) and ii) of Definition 1.14 in [VD-W3].
\snl
Let us now also check that the third set of assumptions in Definition 1.14 of [VD-W3] are fulfilled. For this we need to consider the maps $G_1$ and $G_2$ as defined in Proposition 1.11 of [VD-W3]. Using e.g.\ the first of formula (2.4), one verifies
$$\align (G_1\ot\iota)(\Delta_{13}(a)(1\ot b\ot c))
&=\Delta_{13}(a)(1\ot E)(1\ot b\ot c)\\
&=\Delta_{13}(a)(F_1\ot 1)(1\ot b\ot c)
\endalign$$
for all $a,b,c\in A$. By the fullness of $\Delta$ we conclude that 
$G_1(a\ot b)=(a\ot 1)F_1(1\ot b)$
for all $a,b$. Similarly we will find that 
$G_2(a\ot b)=(a\ot 1)F_2(1\ot b)$
for all $a,b$. Then, the result of Proposition 2.12 provides the necessary condition required in Definition 1.14 for the kernels of the canonical maps.
\snl
This shows that $(A,\Delta)$ is a weak multiplier Hopf algebra. Regularity easily follows from the fact that we can e.g.\ replace $A$ by $A^{\text{op}}$ since we also have the necessary properties for the ranges and the kernels of $T_3$ and $T_4$.  
\hfill$\square$
\einspr

The counit is weakly multiplicative, in the sense of Definition 1.4. This follows, as in the proof of Proposition 4.12 of [VD-W3] from the previous result. This shows that we are in fact in the situation as in the beginning of Section 1. 
\snl
In [VD7], some more aspects of the Larson-Sweedler theorem for weak multiplier Hopf algebras are considered.
\nl
We now complete this section by formulating some properties of the antipode $S$ of the weak multiplier Hopf algebra $(A,\Delta)$.

\iinspr{2.15} Proposition \rm
Assume that $\varphi$ is a left integral. Let $a,b\in A$ and define \newline 
$p=(\iota\ot\varphi)(\Delta(a)(1\ot b))$. Then $S(p)=(\iota\ot\varphi)((1\ot a)\Delta(b))$.

\snl\bf Proof\rm:
We have shown that for any left integral $\varphi$ we have
$$T_1((\iota\ot\iota\ot\varphi)\Delta_{13}(a)\Delta_{23}(b)(1\ot c\ot 1))=E(p\ot c)$$
where $a,b,c\in A$ and where $p=(\iota\ot\varphi)(\Delta(a)(1\ot b))$. When $R_1$ is the generalized inverse of $T_1$ giving the antipode $S$ by the formula 
$$R_1(p\ot c)=\sum_{(p)}p_{(1)}\ot S(p_{(2)})c,$$
we find, using that $R_1T_1$ is given by $R_1T_1(u\ot v)=(u\ot 1)F_1(1\ot v)$ for all $u,v$ that 
$$\align
R_1(p\ot c)
&=R_1T_1((\iota\ot\iota\ot\varphi)\Delta_{13}(a)\Delta_{23}(b)(1\ot c\ot 1)))\\
&=(\iota\ot\iota\ot\varphi)\Delta_{13}(a)(F_1\ot 1)\Delta_{23}(b)(1\ot c\ot 1))\\
&=(\iota\ot\iota\ot\varphi)\Delta_{13}(a)(1\ot E)\Delta_{23}(b)(1\ot c\ot 1))\\
&=(\iota\ot\iota\ot\varphi)\Delta_{13}(a)\Delta_{23}(b)(1\ot c\ot 1)).
\endalign$$
As this is equal to $\sum_{(p)}p_{(1)}\ot S(p_{(2)})c$, we get
$$S(p)c=(\iota\ot\varphi)(1\ot a)\Delta(b))(c\ot 1).$$
We can now cancel $c$ and this completes the proof.
\hfill$\square$
\einspr

In fact, this result can be proven without reference to the previous formulas. Indeed, it is shown to be true for any left integral on a regular weak multiplier Hopf algebra in [VD-W5]. Similarly for other cases.  

\iinspr{2.16} Remark \rm
The above result suggests another approach to the theory. The argument presented in the proof of the previous proposition can be used to {\it define} the antipode. This idea has been used e.g.\ already in [VD5].
\hfill$\square$
\einspr

We finally verify that $S$ coincides with the antipodal maps on the source and target algebras.

\iinspr{2.17} Proposition \rm
The antipode $S$ of $(A,\Delta)$, extended to the multiplier algebra $M(A)$ coincides with the antipodal maps $S_B$ on $B$ and with $S_C$ on $C$.

\snl\bf Proof\rm:
Let $y\in C$. For any pair $a,b$ of elements in $A$ and for any left integral $\varphi$ we have 
$$\align 
(\iota\ot \varphi)((1\ot ay)\Delta(b))
&=S((\iota\ot\varphi)(\Delta(ay)(1\ot b)))\\
&=S((\iota\ot\varphi)(\Delta(a)(y\ot b)))\\
&=S((\iota\ot\varphi)(\Delta(a)(1\ot b))y)\\
&=S(y)S((\iota\ot\varphi)(\Delta(a)(1\ot b)))\\
&=S(y)(\iota\ot\varphi)((1\ot a)\Delta(b)).
\endalign$$
As this is true for all left integrals and all elements $a$, it follows that 
$$(1\ot y)\Delta(b)=(S(y)\ot 1)\Delta(b)$$
for all $b$. This implies that $S(y)=S_C(y)$.
\snl
Similarly, the antipode $S$ of the weak multiplier Hopf algebra $(A,\Delta)$ coincides on $B$ with the antipode map $S_B$ associated with $E$.
\hfill$\square$
\einspr

Also here, we should have a quick look at the involutive case, just as we did at the end of the previous section. However, we postpone this to the next section, where discuss the possible development of a theory in the operator algebra setting.
\nl\nl

\bf 3. Conclusions and further remarks \rm
\nl
This paper has two main sections.
\snl
In the first one, Section 1, we have discussed the arguments needed to show that the (unique) canonical idempotent $E$ in the definition of a weak multiplier bialgebra is a (regular) separability idempotent in the multiplier algebra $M(B\ot C)$ where $B$ and $C$ are the images of the counital maps $\varepsilon_s$ and $\varepsilon_t$. Also these maps are obtained from the axioms of a weak multiplier bialgebra. The result has been proven by B\"ohm, G\'omez-Torrecillas and L\'opez-Centella in [B-G-L]. The existence of a counit $\varepsilon$ is assumed from the very beginning and $\varepsilon$ is required to satisfy two weak multiplicativity axioms. Proofs are not given here as they are found in [B-G-L]. This section is mainly included for completeness and because the result is crucial for the interpretation of the main result of this paper.
\snl
The main result of this paper is the Larson-Sweedler theorem for weak multiplier Hopf algebras. It says that, if a weak multiplier bialgebra has enough integrals, then it has an antipode and it is a (regular) weak multiplier Hopf algebra. This result is obtained by combining the results of the first section with the ones of the second section.
\snl
However, in accordance with the aim of this paper, as explained in the introduction - see also further here - we start in Section 2 with a slightly different setting as the one obtained at the end of Section 1. The main difference lies in the fact that we no longer assume the existence of a counit. We start with a pair of a non-degenerate algebra $A$ and a full and regular coproduct  on $A$. We assume that the canonical idempotent $E$ in $M(A\ot A)$ exists and that it is a separability idempotent. Under these conditions, if there exist enough integrals, we obtain an invertible antipode, as well as a counit satisfying the weak multiplicativity as required in the first section.
\snl
Recently, the theory of separability idempotents for non-degenerate algebras has been considered with less regularity constraints, see [VD4.v2]. Recall that in the first version of this paper, reference [VD4.v1], only regular separability idempotents were considered. This paper on the Larson-Sweedler theorem has been written before the more general non-regular separability idempotents were studied. This raises the question whether or not the Larson-Sweedler theorem as we treat it in this paper, can further be generalized with weaker regularity conditions. This is not at all clear as it seems that the existence of faithful integrals implies regularity in the end. Still, the problem is worth an investigation and a first attempt is found in a paper with reflections on the Larson-Sweedler theorem for (weak) multiplier Hopf algebras, see [VD7]. Remark however that the non-regular case is not relevant when we talk about the involutive setting as then regularity is always automatic.
\nl
The Larson-Sweedler theorem is well-known among algebraists. On the other hand, it plays an implicit role in the operator algebra approach to quantum groups. Consider e.g.\ the definition of a locally compact quantum group in the setting of operator algebras. In the von Neumann algebra framework, the starting point is a pair of a von Neumann algebra with a normal and unital coproduct  and with the assumption that there is a left and a right Haar weight. See e.g.\ [K-V2] and also [VD5]. The existence of a counit, as well as of an antipode, is not part of the starting assumptions. This is just like in the case of the Larson-Sweedler theorem, as we have treated it in Section 2 of this paper.
\snl
Indeed, this is precisely the situation that we start with in Section 2, but now for weak multiplier Hopf algebras. And the challenge is clear. The aim is to develop the theory as it now exists for weak multiplier Hopf algebras in the operator algebra framework, either in the C$^*$-algebra setting or in the von Neumann algebra context. What is done in Section 2 of this paper should be a possible source of inspiration for such a project (as it was in the case of the development of locally compact quantum groups).
\snl
Of course, first one has to see what happens with the results in Section 2 in the case of a $^*$-algebra. In Appendix B, where we recall the notion of a separability idempotent, we briefly mention this case while in the original paper ([VD4.v1]), this case got special attention. However not in Section 1, nor in Section 2, we have really investigated the involutive case yet. Therefore let us have a look at both the starting point and the final result in Section 2, for the involutive case and give some comments.
\snl
So consider a non-degenerate idempotent $^*$-algebra $A$ and a full coproduct $\Delta$ on $A$. Assume that $\Delta$ is a $^*$-map. It will follow that the canonical idempotent $E$ is self-adjoint because it is the smallest idempotent in $M(A\ot A)$ satisfying 
$E\Delta(a)=\Delta(a)E=\Delta(a)$ 
for all $a\in A$. Remark also that the equality 
$$(\Delta\ot\iota)E=(1\ot E)(E\ot 1)$$
will imply 
$$(1\ot E)(E\ot 1)=(E\ot 1)(1\ot E),$$
again by taking adjoints.
\snl
The separability idempotent $E$ in $M(B\ot C)$, being self-adjoint, will imply that $B$ and $C$ are $^*$-subalgebras of $M(A)$. Then we know from the theory of separability idempotents, as developed in [VD4.v1], that $B$ and $C$ are direct sums of full matrix algebras. Moreover, the functionals $\varphi_B$ on $B$ and $\varphi_C$ on $C$ are positive. Furthermore, the pair of the algebra $B$ and the faithful positive linear functional $\varphi_B$ on $B$ determine the separability idempotent completely. Recall from [VD4.v2] that in the involutive case, the separability idempotent is automatically regular. 
\snl
We will need the extra assumption that there are enough {\it positive} integrals on $A$. As these integrals are not assumed to be faithful on their own, there may not exist modular automorphisms. However, they still have some 'restricted modular behavior' in the sense that e.g.\ $\varphi(ya)=\varphi(a\sigma_C(y))$ when $\varphi$ is left invariant, $a$ is in $A$ and $y$ in $C$. Here $\sigma_C$ is the modular automorphism of $\varphi_C$ on $C$. See a remark after the proof of Proposition 2.7. It might happen that, in the involutive case, if enough integrals exists, there are automatically enough positive integrals. This however is not clear.
\snl
All this suggests a possible set of axioms for the generalization of the concept of a weak multiplier Hopf $^*$-algebra to the operator algebra setting. In the framework of von Neumann algebras, the starting point will be a von Neumann algebra $M$ and a coproduct $\Delta$ on $M$ satisfying the obvious conditions (without being unital). There is also the subalgebra $N$ of $M$ with a faithful normal semi-finite weight $\nu$. The pair $(N,\nu)$ has to give rise to the canonical idempotent $E$, now a self-adjoint projection in the von Neumann algebraic tensor product of $N\ot L$ where $L$ is a von Neumann subalgebra of $M$, anti-isomorphic with $N$. Finally, there will be the requirement of having enough normal semi-finite (not necessarily faithful) weights on $M$ that are left invariant, as well as having enough such weights that are right invariant. They will need to have some restricted modular behavior in the sense given by the modular automorphisms of $\nu$ on $N$. This project is carried out further in [K-VD2] and [K-VD3], see also [K-VD1]. 
\snl
Finally we should also say something about the relation of all of this with the existing theory of measured quantum groupoids as developed e.g.\ by Enock and Lesieur in [L] and [E]. The difference of their work and the theory as outlined above lies in the existence of the canonical idempotent $E$ for the pair $(N,\nu)$. This is of the same nature as the difference between multiplier Hopf algebroids (as developed in [T-VD1]) and the weak multiplier Hopf algebras as developed in [VD-W3]. One of the main conditions for a multiplier Hopf algebroid to have an underlying weak multiplier Hopf algebra is that the base algebra is separable Frobenius, see [T-VD2]. And so it is expected that there should be some equivalent condition in the case of operator algebras.
\nl\nl

\bf Appendix A. A dictionary\rm
\nl
In analysis, in particular in functional analysis and the theory of operator algebras, several attempts to generalize Pontryagin's duality for abelian locally compact groups to the non-abelian case have led to the interest of Hopf algebraic structures. More than a decade ago, this resulted in a theory of locally compact quantum groups, see e.g.\ [K-V1] and [K-V2]. At the same time, in fact even earlier, the algebraists were studying Hopf algebras and more recently quantum groups in a purely algebraic setting. Unfortunately, both fields use different terminologies, different notations and sometimes even a very different way of thinking and arguing. 
\snl
A similar phenomenon occurred with the introduction of the theory of quantum groupoids. In pure algebra, there is the literature that started with the study of weak Hopf algebras by G.\ B\"ohm, F.\ Nill \& K.\ Szlach\'anyi [B-N-S]. Also the operator algebraists were interested because of the relation with the study of subfactors (and partly also as a generalization of the theory of locally compact quantum groups to locally compact quantum groupoids). And although in the very beginning, the original work was entitled: {\it Weak Hopf algebras I. Integral theory and C$^*$-structure}, soon after, the algebraists developed the theory of weak Hopf algebras independently of the study of these objects in the operator algebras setting.
\snl 
We feel that the theory of weak multiplier Hopf algebras, and in particular the theory of weak multiplier Hopf $^*$-algebras with positive integrals, can be considered to be a bridge between the two approaches. Such was already the case with multiplier Hopf algebras, more precisely with the theory of multiplier Hopf $^*$-algebras with positive integrals, sometimes called algebraic quantum groups. An important aspect of this last theory is the possibility to construct duals in very much the same way as is done in the original Pontryagin's duality for abelian locally compact groups. Moreover, this duality extends the duality of finite-dimensional Hopf algebras to a more general, infinite-dimensional case.
\snl
The aim of this appendix is to provide a small dictionary as far as the topic of this paper is concerned, namely about weak multiplier bialgebras and weak multiplier Hopf algebras. First we look at the more general aspects.
\nl
\it Differences and common practices \rm
\nl
First, in analysis and certainly in operator algebras, most if not all algebras that play a role, are {\it algebras over the field of complex numbers}, whereas in algebra, not only algebras over other fields are considered, but often people just work with rings. This implies that several results about e.g.\ tensor products, obvious for analysts, are no longer valid.
\snl
In analysis one is often also interested in the case where the underlying algebras are {\it $^*$-algebras} and mostly where they are {\it operator algebras}. This means that they have a $^*$-representation by say bounded operators on a Hilbert space. This will imply e.g.\ that $a^*a=0$ can only happen if $a=0$.
\snl
The notion of a multiplier algebra of a non-unital, but non-degenerate algebra is something familiar for people working with operator algebras. This is not so for algebraists. For algebraists, a multiplier of an algebra would be considered as a pair $(\lambda,\rho)$ of module maps from  $A$ to itself, satisfying the compatibility $\rho(a)b=a\lambda(b)$ for all $a,b\in A$. This is also true for operator algebraists. But they will rather look at multipliers as sitting in the multiplier algebra $M(A)$, characterized as the largest algebra with identity in which $A$ sits as an essential ideal. Then $A$ is seen as a subalgebra of $M(A)$. This is however only possible if the product in the algebra is non-degenerate and so usually, in this setting, only such algebras are considered.
\snl
For a {\it linear functional} $\omega$ on an algebra $A$, there is the notion of {\it faithfulness} as used in analysis. It means that an element $a\in A$ has to be $0$ if $\omega(ab)=0$ for all $b\in A$ and also if $\omega(ba)=0$ for all $b\in A$. In the case of a $^*$-algebra, and when the functional is positive, this is equivalent with the property that $a=0$ if $\omega(a^*a)=0$. The existence of a faithful functional on an algebra implies that the product is non-degenerate. The existence of a positive faithful functional on a $^*$-algebra is related with the property that the algebra is an operator algebra, although this is not sufficient to have a faithful $^*$-representation by {\it bounded} operators on a Hilbert space. Fortunately, in the situations we consider in Hopf algebra theory, it seems to be so that e.g.\ positivity of the integrals give rise to bounded operators on the associated Hilbert spaces.
\snl
For a finite-dimensional algebra the existence of a faithful linear functional implies not only that the algebra is unital, but also that it is Frobenius. Also the converse is true here. The Frobenius property is not familiar to operator algebraists because there, often algebras are nice $^*$-algebras and the property is always satisfied.
\snl
If $\omega$ is a faithful linear functional on a finite-dimensional algebra $A$, there automatically exists an automorphism $\sigma$ satisfying $\omega(ab)=\omega(b\sigma(a))$ for all $a,b\in A$. This result is no longer true for infinite-dimensional algebras. Nevertheless, in the situation studied here, the integrals turn out to have this property, also in the infinite-dimensional case. In the field of functional analysis, this automorphism is called the {\it modular automorphism}. The terminology finds its origin in the theory of non-unimodular locally compact groups and by extension of this, in the modular theory of faithful positive linear functionals on a von Neumann algebra (see e.g.\ [T]). In algebra, this automorphism (or rather its inverse) is known as the {\it Nakayama automorphism} (see e.g.\ [Na]).
\snl
For the identity map on the space $A$, we intend to use $\iota_A$ or $\text{id}_A$ and often we even drop the index because most of the time it is clear from the context on which space the identity map $\iota$ acts. In algebra, people often use the same symbol for the space and for the identity map on this space. For a coproduct $\Delta:A\to A\ot A$, whereas we use e.g.\ $\Delta\ot\iota$ for the map from $A\ot A$ to $A\ot A\ot A$, obtained by applying $\Delta$ to the first factor and the identity map on the second factor, in the algebra literature, it is common to write this map as $\Delta\ot A$. This looks quite strange to functional analysts. 
\snl
For algebraists, it is quite natural to consider the dual notion, namely that of a {\it coalgebra}. It is  a vector space $C$ with a coproduct $\Delta:C\to C\ot C$ that is coassociative and has a counit. In many situations however this notion of a coproduct is {\it too restrictive} and one needs a bigger space for the range of the coproduct. If e.g.\ $A$ is a non-degenerate algebra, a  natural concept of a coproduct on $A$ is a homomorphism $\Delta:A\to M(A\ot A)$ and in order to be able to formulate coassociativity, there is the need of some extra assumptions about the range of $\Delta$ in $M(A\ot A)$ (see the discussion in the beginning of Section 1). There is also the notion of regularity and fullness of a coproduct that is completely irrelevant if the coproduct maps into the tensor product and if there is a counit.
\snl
There seems to be no way to extend the notion of a coproduct simply on a vector space that fits our needs, without more structure.
\nl
\it More specific comments about weak multiplier bialgebras \rm 
\nl
As homomorphisms like $\Delta$ go from one algebra into the multiplier algebra of another one, compositions of such maps are not immediately possible. One has to find a way to extend such homomorphisms to the multiplier algebra. There is a general theory about this for {\it non-degenerate} homomorphisms $\gamma$ from an algebra $A$ to the multiplier algebra $M(B)$ of an algebra $B$. Recall that $\gamma$ is called non-degenerate if $\gamma(A)B=B\gamma(A)=B$. In that case, there exists a unique extension $\overline\gamma$ which is a unital homomorphism from $M(A)$ to $M(B)$. In some sense, the non-degenerate homomorphisms from $A$ to $M(B)$ are the natural generalizations to non-unital non-degenerate algebras of unital homomorphisms in the case of unital algebras. However, this is only true for idempotent algebras. It seems to be more subtle when the algebras are not idempotent. See some remarks about this problem in [VD8].
\snl
In the theory of weak multiplier Hopf algebras, the coproduct (in general) is not a non-degenerate homomorphism from $A$ to $M(A\ot A)$ and so the above extension procedure does not work. Fortunately, there is the canonical idempotent $E$ with its properties that allows to find such extensions. This is explained e.g.\ in an appendix of [VD-W2]. See also the remark in the beginning of Section 1 of this paper.
\snl
It is also a common practice in the field of operator algebras, to use the same symbols for these extensions as we do in this paper as well.
\snl
In Definition 1.2 we use the term {\it weakly non-degenerate} for the coproduct, whereas in the original work on weak Hopf algebras [B-N-S] this is called the {\it weak multiplicativity of the unit}. The latter is quite natural for algebras that have an identity whereas the former is perhaps more natural if the algebra has no unit because then it is indeed a weaker form of non-degeneracy and it is used in the same way to extend homomorphisms. See also some comments in Section 1 on the use of this terminology.
\snl
Now we consider the {\it counital maps}. In the original paper by B\"ohm, Nill \& Szlach\'anyi, in the case of a weak Hopf algebra, the counital maps are denoted and defined as follows in terms of the antipode $S$. The Sweedler notation $\sum_{(a)}a_{(1)}\ot a_{(2)}$ is used for $\Delta(a)$. 

\inspr{A.1} Notation \rm Let $(A,\Delta)$ be a weak Hopf algebra with antipode $S$. Then we have the maps
$$\align
\sqcap^R(a)&=\sum_{(a)}S(a_{(1)})a_{(2)}
\qquad\qquad\qquad\qquad
\overline\sqcap^R(a)=\sum_{(a)}a_{(2)}S^{-1}(a_{(1)})\\
\sqcap^L(a)&=\sum_{(a)}a_{(1)}S(a_{(2)})
\qquad\qquad\qquad\qquad
\overline\sqcap^L(a)=\sum_{(a)}S^{-1}(a_{(2)})a_{(1)}
\endalign$$
\hfill $\square$
\einspr

These counital maps can also be defined for a weak multiplier bialgebra, in terms of the counit, without reference to the (non-existing) antipode. We have the formulas in Proposition 1.7 and we recall them here. We have, for all $a\in A$,
$$\align 
\varepsilon_s(a)&=(\iota\ot\varepsilon)((1\ot a)E)
\qquad\qquad\qquad
\varepsilon_s'(a)=(\iota\ot\varepsilon)(E(1\ot a))\\
\varepsilon_t(a)&=(\varepsilon\ot\iota)(E(a\ot 1))
\qquad\qquad\qquad
\varepsilon_t'(a)=(\varepsilon\ot\iota)((a\ot 1)E).
\endalign$$
Now, here are the connections with the case of a weak (multiplier) bialgebra and the formulas as found e.g.\ in Section 3 of [B-G-L]. For $a\in A$ we have
$$\align 
\varepsilon_s(a)&=\sqcap^R(a)
\qquad\qquad\qquad\qquad\qquad
\varepsilon_s'(a)=\overline\sqcap^R(a)\\
\varepsilon_t(a)&=\sqcap^L(a)
\qquad\qquad\qquad\qquad\qquad
\varepsilon_t'(a)=\overline\sqcap^L(a).
\endalign$$
The notations we use for the counital maps, are also used in e.g.\ [N-V2].
\nl
In the next appendix, we will say more about the common conventions and differences, also differences in notations, for the notion of {\it separability} and {\it separability idempotents}.
\nl\nl

\bf Appendix B. Separability and Frobenius algebras \rm
\nl
In this appendix we will first recall the notion of a {\it (regular) separability idempotent} as it is studied in [VD4.v1]. As we mentioned already, in the more recent version [VD4.v2], also non-regular separability idempotents are considered. However, in this note, only the regular case is studied and therefore, it is sufficient to recall here only the material that is already present in the first version of the paper on separability for multiplier algebras.
\snl
We not only recall the main properties, we will also prove some new properties we need in this paper. Throughout we will give some comments, refer to the existing literature about separable algebras and {\it continue with the dictionary} of the previous appendix.
\snl
We start with two non-degenerate algebras $B$ and $C$ and an idempotent $E$ in $M(B\ot C)$. We require that 
$$\align
&E(1\ot c)\in B\ot C \qquad\quad\text{and}\qquad\quad (b\ot 1)E\in B\ot C \tag"(B.1)"\\
&(1\ot c)E\in B\ot C  \qquad\quad\text{and}\qquad\quad E(b\ot 1)\in B\ot C \tag"(B.2)"
\endalign $$
for all $b\in B$ and $c\in C$. It is further assumed that the left and the right leg of $E$ are respectively all of $B$ and $C$. This means that $E$ is full in the terminology of [VD4.v1]. 
\snl
Remark that in the setting of Section 2, we have the element $E$ sitting in $M(A\ot A)$ and we also have that 
$$\align
&E(1\ot a)\in B\ot A \qquad\quad\text{and}\qquad\quad (a\ot 1)E\in A\ot C \\
&(1\ot a)E\in B\ot A \qquad\quad\text{and}\qquad\quad E(a\ot 1)\in A\ot C 
\endalign $$
for all $a\in A$. See the proof of Proposition 2.1.
\snl 
In what follows, we will systematically use letters $b,b'$ for elements in $B$ and letters $c,c'$ for elements in $C$. We can do this here because further, we do not consider the bigger algebra $A$.
\snl
We now recall the following definition from [VD4.v1].

\inspr{B.1} Definition \rm
We call $E$ a {\it regular separability idempotent} if also
$$E(B\ot 1)=E(1\ot C) 
\qquad\quad\text{and}\quad\qquad
(B\ot 1)E=(1\ot C)E. \tag"(B.3)"$$
\vskip -0.7 cm	
\hfill$\square$\einspr

Remark once more that in [VD4.v1] all separability idempotents are assumed to be regular in the sense of the generalization studied in the second version [VD4.v2].
\snl
One can show that this condition implies that the algebras have {\it local units} (see Proposition 1.9 in [VD4.v1]). In particular they are idempotent algebras. 
\snl
Also remark that the conditions (B.1) and (B.2) follow from (B.3) and the fact that $E$ is a multiplier of $B\ot C$ (when $B$ and $C$ are already known to be idempotent). 
\snl
We have the existence of the {\it antipodal maps} $S_B$ and $S_C$:

\inspr{B.2} Proposition \rm
Assume that $E$ is a regular separability idempotent in $M(B\ot C)$. Then there are anti-isomorphisms $S_B:B\to C$ and $S_C:C\to B$ given by 
$$E(b\ot 1)=E(1\ot S_B(b)) \qquad\quad\text{and}\qquad\quad (1\ot c)E=(S_C(c)\ot 1)E$$
for all $b\in B$ and $c\in C$.
\hfill$\square$\einspr

See Proposition 1.7 in [VD4.v1]. 
\snl
Again, from the existence of these anti-isomorphisms, the previous condition (B.3) follows. We see that an idempotent $E\in M(B\ot C)$ will be a regular separability idempotent if the algebras $B$ and $C$ are idempotent, if the anti-isomorphisms $S_B$ and $S_C$ exist and if $E$ is full. 
\snl
As $B$ is {\it anti-isomorphic} with $C$, we could have taken $C=B^{\text{op}}$ from the very beginning (using e.g.\ $S_C$ - see remark B.4.i below). One of the reasons for not doing this is the fact that we are also interested in the involutive case where in general $S_C:C\to B$ is not a $^*$-anti-isomorphism. See also some more comments on the involutive case at the end of this appendix.
\snl
It is a fairly easy consequence that the following holds (see Proposition 1.8 in [VD4.v1]).

\inspr{B.3} Proposition \rm
If $E$ is a regular separability idempotent in $M(B\ot C)$, then 
$$m_C(S_B\ot\iota)(E(1\ot c))=c
\qquad\qquad\text{and}\qquad\qquad
m_B(\iota\ot S_C)((b\ot 1)E)=b$$ 
for all $c\in C$ and all $b\in B$ where $m_B$ and $m_C$ are used to denote the multiplication maps from $B\ot B\to B$ and $C\ot C\to C$ respectively. 
\hfill$\square$\einspr

We think of the formulas above as
$$S_B(E_{(1)})E_{(2)}=1 \qquad\quad\text{and}\qquad\quad E_{(1)}S_C(E_{(2)})=1$$
with the Sweedler type notation $E=E_{(1)}\ot E_{(2)}$.   
\snl

One can also show that $(S_B\ot S_C)E=\zeta E$ where $\zeta$ is the flip from $B\ot C$ to $C\ot B$, extended to the multiplier algebra (Proposition 1.13 in [VD4.v1]).
\nl
Before we continue, let us compare this notion with the notion of separability as it appears elsewhere in the literature (see e.g.\ [DM-I]). 

\inspr{B.4} Remark \rm
i) Consider the case where $B$ and $C$ are finite-dimensional algebras with an identity, so that in particular $E\in B\ot C$. The anti-isomorphism $S_C$ can be used to identify $C$ with $B^{\text{op}}$. Define $F=(\iota\ot S_C)E$. Then $F\in B\ot B^{\text{op}}$ and we will have $mF=1$ as well as $(b\ot 1)F=F(1\ot b)$ for all $b\in B$. This means that $F$ is a separability idempotent in the sense of the literature. 
\snl
ii) Moreover, in our case, the two legs of $F$ are all of the algebra $B$. This is not necessarily true for any separability element. Consider e.g.\ $B=M_n(\Bbb C)$, the algebra of $n\times n$ complex matrices, and let 
$$F=\sum_i e_{i1}\ot e_{1i}$$
where $e_{ij}$ are matrix elements. Then $F$ is a separability element in the sense of i), but the legs do not give all of the algebra. Therefore, it will not fit into the concept as we studied it in [VD4.v1] (and as we use it here).
\hfill$\square$\einspr

The separability idempotents as we use them here are all of \it Frobenius type\rm.
\snl
Now we again consider the general setting and any regular separability idempotent $E$ as in Definition B.1 above.

\inspr{B.5} Proposition \rm There exist unique linear functionals $\varphi_B$  on $B$ and $\varphi_C$ on $C$ so that 
$$(\varphi_B\ot\iota)E=1
\qquad\qquad\text{and}\qquad\qquad
(\iota\ot\varphi_C)E=1$$
in $M(C)$ and $M(B)$ respectively. These functionals are faithful.
\hfill$\square$\einspr

It follows from the formula $(S_B\ot S_C)E=\zeta E$ and the uniqueness that $\varphi_B=\varphi_C\circ S_B$ and $\varphi_C=\varphi_B\circ S_C$. 
\snl
If we now consider again the finite-dimensional case and $F=(\iota\ot S_C)E$ as in the remark, we see that 
$$(\varphi_B\ot\iota)F=1
\qquad\quad\text{and}
\qquad\quad
(\iota\ot\varphi_B)F=1.$$

Because it can be shown that $B$ is made into a coalgebra with the coproduct $\Delta_B$ on $B$ defined as 
$$\Delta_B(b)=(b\ot 1)F=F(1\ot b),$$
in the literature,  $\varphi_B$ is called a counit. Indeed, we clearly have 
$$(\varphi_B\ot\iota)\Delta_B(b)=b
\qquad\quad\text{and}
\qquad\quad
(\iota\ot\varphi_B)\Delta_B(b)=b$$
with this definition of $\Delta_B$.
\snl
In [VD4.v1], we have called $\varphi_B$  and $\varphi_C$ the {\it integrals} on $B$ and $C$ respectively. This seems to be completely contradictory. But in fact, there are arguments to support both points of view, see e.g.\ [VD4.v1] where this is explained. In the more recent version of this paper [VD4.v2] we no longer use this terminology, but call the linear functionals {\it distinguished} linear functionals.
\nl
Again in the finite-dimensional case, the existence of a faithful linear functional is equivalent with the algebra being Frobenius. So if a separability idempotent exists in our terminology, this implies that the algebra is separable and Frobenius. However, the separability and the Frobenius property are connected and so we should really say that we have a {\it separable Frobenius algebra}. See e.g.\ [Sz], [Sc] or [K-S].
\nl
The distinguished linear functionals  on $B$ and $C$ satisfy the {\it weak K.M.S.-property}. Indeed we have 
$$\varphi_B(bb')=\varphi_B(b'\sigma_B(b))
\qquad\quad\text{and}\qquad\quad
\varphi_C(cc')=\varphi_C(c'\sigma_C(c))$$
where $\sigma_B(b)=S_CS_B(b)$ and $\sigma_C(c)=S_B^{-1}S_C^{-1}(c)$ (see Proposition 2.3 in [VD4.v1]). Such automorphisms automatically exist for any faithful functional on a finite-dimensional algebra. They are then called the Nakayama automorphisms. More precisely, the Nakayama automorphism $\theta$ for the linear functional $\varphi_B$ on $B$ is the inverse of $\sigma_B$.  See also some of our remarks about this topic in the previous appendix.
\snl
As it turns out, these modular automorphisms also exist for the integrals on a general separability idempotent. We tend to call these automorphisms the {\it modular automorphisms}. The terminology is motivated by the modular function relating the left and right Haar measure on a non-unimodular locally compact group (and furthermore also by the modular theory in operator algebras).
\nl
Finally, we consider the multipliers $F_i$ used in Propositions 2.6 and 2.7 of this paper. We also give a proof of Proposition 2.6.

\inspr{B.6} Proposition \rm
There exist elements $F_1$ and $F_3$ in $M(B\ot B)$ and  $F_2$ and $F_4$ in $M(C\ot C)$ characterized by 
$$\align F_1&=(\iota\ot S_C)E  \qquad\quad\text{and}\qquad\quad F_3=(\iota\ot S_B^{-1})E\\
	F_2&=(S_B\ot \iota)E \,\qquad\quad\text{and}\qquad\quad F_4=(S_C^{-1}\ot \iota)E.
\endalign$$

\snl\bf Proof\rm:
Because we assume that $(b\ot 1)E$ and $E(b\ot 1)$ belong to $B\ot C$ for all $b$, we can define the element $F_1$ in $M(B\ot B)$ by
$$(b\ot 1)F_1=(\iota\ot S_C)((b\ot 1)E)
\qquad\quad\text{and}\qquad\quad
F_1(b\ot 1)=(\iota\ot S_C)(E(b\ot 1).$$
Similarly for the other cases.
\hfill$\square$\einspr

We see that actually we already have that $(b\ot 1)F_1$ and $F_1(b\ot 1)$ belong to $B\ot B$ for all $b\in B$. Also  $(1\ot b)F_1$ and $F_1(1\ot b)$ belong to $B\ot B$ for all $b\in B$. Similarly for the other elements $F_2$, $F_3$ and $F_4$.
\snl
Now we can easily prove the result of Proposition 2.6 of Section 2 of this paper.

\inspr{B.7} Proposition \rm  The elements $F_i$ satisfy 
$$\align E_{13}(F_1\ot 1) &= E_{13}(1\ot E) 
	\quad\text{ and }\quad 
		(F_3\ot 1)E_{13}=(1\ot E)E_{13} \tag"(B.4)"\\
	(1\ot F_2)E_{13} &=(E\ot 1)E_{13}
	\quad\text{ and }\quad
	 	E_{13}(1\ot F_4)=E_{13}(E\ot 1).\tag"(B.5)"
\endalign$$

\bf\snl Proof\rm:
Take any element $b$ in $B$. If we use the Sweedler type notation
$E_{(1)}b\ot E_{(2)}$ for $E(b\ot 1)$ we find
$$E_{13}(1\ot E)(1\ot b\ot 1)=E_{13}(S_B^{-1}(E_{(2)})\ot E_{(1)}b\ot 1).$$
Because $(S_B\ot S_C)E=\zeta E$ we have 
$$S_B^{-1}(E_{(2})\ot E_{(1)}b=E_{(1)}\ot S_C(E_{(2)})b$$
and this proves the first formula.
\snl
The other formulas are proven in a completely similar way.
\hfill$\square$
\einspr

Observe that the expressions in (B.4) are in $M(B\ot B\ot C)$ whereas the elements in the formulas (B.5) are in $M(B\ot C\ot C)$. We also have that the elements $F_i$ are determined by these formulas (because $E$ is full). 
\nl
We finish this appendix with some more information about the {\it involutive case} (cf. Section 3 in [VD4.v1]). In this setting, it is assumed that $B$ and $C$ are $^*$-algebras and that $E$ is self-adjoint. It follows that
$$S_C(S_B(b^*))^*=b
\quad\quad\text{and}\quad\quad 
S_B(S_C(c)^*)^*=c$$
for all $b\in B$ and $c\in C$. The linear functionals $\varphi_B$ and $\varphi_C$ are positive. We also see that the antipodal maps $S_B$ and $S_C$ are $^*$-maps if and only if the modular automorphisms $\sigma_B$ and $\sigma_C$ are trivial, that is if $\varphi_B$ and $\varphi_C$ are traces on $B$ and $C$ respectively. Because this is not true in general, from this point of view, it is not appropriate to identify $C$ with $B^{\text{op}}$ because they are not identified as involutive algebras.
\nl\nl

\bf References \rm
\nl
{[\bf Abe]} E.\ Abe: {\it Hopf algebras}. \rm Cambridge University Press (1977).
\snl
{[\bf Abr]} L.\ Abrams: {\it Modules, comodules and cotensor products over Frobenius algebras}. J.\ Algebra 219 (1999), 201-213.
\snl
{[\bf B-N-S]} G.\ B\"ohm, F.\ Nill \& K.\ Szlach\'anyi: {\it Weak Hopf algebras I. Integral theory and C$^*$-structure}. J.\ Algebra 221 (1999), 385-438. 
\snl
{[\bf B-G-L]} G.\ B\"ohm, J.\ G\'omez-Torecillas and E.\ L\'opez-Centella: {\it Weak multiplier bialgebras}.  Trans. Amer. Math. Soc. 367 (2015), no. 12, 8681-872.  See also arXiv: 1306.1466 [math.QA]. 
\snl
{[\bf B-S]} G.\ B\"ohm  \& K.\ Szlach\'anyi: {\it Weak Hopf algebras II. Representation theory, dimensions and the Markov trace}. J.\ Algebra 233 (2000), 156-212. 
\snl
{[\bf DM-I]} F.\ De Meyer \& E. Ingraham: {\it Separable algebras over a commutative rings}. Lecture Notes in Mathematics 181 (1971).
\snl
{[\bf E]} M.\ Enock: {\it Measured quantum groupoids with a central basis}. J.\ Operator Theory, 66 (2011), 3-58.
\snl
{[\bf J-V]} K.\ Janssen \& J.\ Vercruysse: {\it Multiplier bi- and Hopf algebras}.  J.\ Algebra Appl.\ 9 (2), (2010), 275-303. 
\snl
{[\bf K-S]} L.\ Kadison \& K.\ Szlach\'anyi: {\it Bialgebroid actions on depth two extensions and duality}. Adv.\ in Math.\  179 (2003),  75-121.
\snl
{[\bf K-V1]} J.\ Kustermans \& S.\ Vaes: \it Locally compact quantum groups. \rm Ann.\ Sci.\ \'Ec.\ Norm.\ Sup.\   33 (2000), 837--934.
\snl
{[\bf K-V2]} J.\ Kustermans \& S.\ Vaes: {\it Locally compact quantum groups in the von Neumann algebraic setting}. Math.\ Scand.\ 92 (2003), 68-92.
\snl
{[\bf K-VD1]} B.-J.\ Kahng \& A.\ Van Daele: {\it Seperability idempotents in C$^*$-algebras }. Preprint Canisius College Buffalo (USA) and University of Leuven (Belgium). 
\snl
{[\bf K-VD2]} B.-J.\ Kahng \& A.\ Van Daele: {\it A class of C$^*$-algebraic locally compact quantum groupoids I }. Preprint Canisius College Buffalo (USA) and University of Leuven (Belgium). 
\snl
{[\bf K-VD3]} B.-J.\ Kahng \& A.\ Van Daele: {\it A class of C$^*$-algebraic locally compact quantum groupoids II }. Preprint Canisius College Buffalo (USA) and University of Leuven (Belgium).
\snl
{[\bf L-S]} R.G.\ Larson \& M.E.\ Sweedler: {\it An associative orthogonal bilinear form for Hopf algebras}. Amer.\ J.\ Math.\ 91 (1969), 75-93.
\snl
{[\bf L]} F.\ Lesieur: {\it Measured quantum groupoids}. M\'emoirs de la Soci\'et\'e Mathematique de France 109 (2007), 1-122.
\snl
{[\bf Na]} T.\ Nakayama: {\it On Frobeniusean algebras II}. Ann.\ of Math.\ 42, (1941), 1-21.
\snl
{[\bf Ni]} D.\ Nikshych: {\it On the structure of weak Hopf algebras}. Adv.\ Math.\ 170 (2002), 257-286.
\snl
{[\bf N-V1]} D.\ Nikshych \& L.\ Vainerman: {\it Algebraic versions of a finite dimensional quantum groupoid}. Lecture Notes in Pure and Applied Mathematics 209 (2000), 189-221. 
\snl
{[\bf N-V2]} D.\ Nikshych \& L.\ Vainerman: {\it Finite quantum groupoids and their applications}. In {\it New Directions in Hopf algebras}. MSRI Publications, Vol.\ 43 (2002), 211-262.
\snl
{[\bf Sw]} M.\ Sweedler: {\it Hopf algebras}. Benjamin, New-York (1969).
\snl
{[\bf Sc]} P.\ Schauenburg: {\it Weak Hopf algebras and quantum groupoids}. Noncommutative geometry and quantum groups (Warsaw, 2001),
 Polish Acad. Sci., Warsaw, (2003), 171-188.
\snl
{[\bf St]} R.\ Street: {\it Frobenius monads and pseudomonoids}. J.\ Math.\ Phys.\ 45 (2004), 3930-3948.
\snl
{[\bf Sz]} K.\ Szlach\'anyi: {\it Finite quantum groupoids and inclusions of finite type}. Fields Inst.\ Comm.\ 30 (2001), 393-407.
\snl
{[\bf T]} M.\ Takesaki: {\it Theory of Operator Algebras II}. Springer, Berlin (2003).
\snl
{[\bf T-VD1]} T.\ Timmermann \& A.\ Van Daele: {\it Regular multiplier Hopf algebroids. Basic theory and examples}. Preprint University of M\"unster (Germany) and University of Leuven (Belgium). See arXiv: 1307.0769 [math.QA]
\snl
{[\bf T-VD2]} T.\ Timmermann \& A.\ Van Daele: {\it Multiplier Hopf algebroids arising from weak multiplier Hopf algebras}. Banach Center Publications, 106 (2015), 73-110. See also arXiv: 1406.3509 [math.RA] 
\snl
{[\bf Va]} L.\ Vainerman (editor): {\it Locally compact quantum groups and groupoids}. IRMA Lectures in Mathematics and Theoretical Physics 2, Proceedings of a meeting in Strasbourg, de Gruyter (2003).
\snl
{[\bf VD1]} A.\ Van Daele: {\it Multiplier Hopf algebras}. Trans.\ Am.\ Math.\ Soc.\  342(2) (1994), 917-932.
\snl
{[\bf VD2]} A.\ Van Daele: {\it An algebraic framework for group duality}. Adv.\ in Math.\ 140 (1998), 323-366.
\snl
{[\bf VD3]} A.\ Van Daele: {\it Tools for working with multiplier Hopf algebras}.   ASJE (The Arabian Journal for Science and Engineering) C - Theme-Issue 33 (2008), 505--528.  See also arXiv:0806.2089 [math.QA].
\snl
{[\bf VD4.v1]} A.\ Van Daele: {\it Separability idempotents and multiplier algebras}. Preprint University of Leuven (Belgium) (2013). See arXiv:1301.4398v1 [math.RA]
\snl
{[\bf VD4.v2]} A.\ Van Daele: {\it Separability idempotents and multiplier algebras}. Preprint University of Leuven (Belgium) (2015). See arXiv:1301.4398v2 [math.RA]
\snl
{[\bf VD5]} A.\ Van Daele: {\it Locally compact quantum groups. A von Neumann algebra approach}. Preprint University of Leuven (Belgium) (2013). SIGMA 10 (2014), 082, 41 pages. See also arXiv:0602.212v2 [math.OA]. 
\snl 
{[\bf VD6]} A.\ Van Daele: {\it The Sweedler notation for (weak) multiplier Hopf algebras}. Prepint University of Leuven (Belgium) (2016). In preparation.
\snl
{[\bf VD7]} A.\ Van Daele: {\it Reflections on the Larson-Sweedler theorem for (weak) multiplier Hopf algebras}. Prepint University of Leuven (Belgium) (2016). In preparation.
\snl
{[\bf VD8]} A.\ Van Daele: {\it A note on coassociativity}. Prepint University of Leuven (Belgium) (2016). In preparation.
\snl
{[\bf VD-W1]} A.\ Van Daele \& S.\ Wang: {\it The Larson-Sweedler theorem for multiplier Hopf algebras}. J.\ of Alg.\ 296 (2006), 75-95.
\snl
{[\bf VD-W2]} A.\ Van Daele \& S.\ Wang: {\it Weak multiplier Hopf algebras. Preliminaries, motivation and basic examples}. Operator Algebras and Quantum Groups. Banach Center Publications 98 (2012), 367-415. See also arXiv:1210.3954v1 [math.RA]
\snl 
{[\bf VD-W3]} A.\ Van Daele \& S.\ Wang: {\it Weak multiplier Hopf algebras I. The main theory}. Journal f\"ur die reine und angewandte Mathematik (Crelles Journal) 705 (2015), 155-209, ISSN (Online) 1435-5345, ISSN (Print) 0075-4102, DOI: 10.1515/crelle-2013-0053, July 2013. See also arXiv:1210.4395v1 [math.RA]
\snl
{[\bf VD-W4.v1]} A.\ Van Daele \& S.\ Wang: {\it Weak multiplier Hopf algebras II. The source and target algebras}. Preprint University of Leuven (Belgium) and Southeast University of Nanjing (China) (2014).  See arXiv:1403.7906v1 [math.RA] 
\snl
{[\bf VD-W4.v2]} A.\ Van Daele \& S.\ Wang: {\it Weak multiplier Hopf algebras II. The source and target algebras}. Preprint University of Leuven (Belgium) and Southeast University of Nanjing (China) (2015).  See arXiv:1403.7906v2 [math.RA] 
\snl
{[\bf VD-W5]} A.\ Van Daele \& S.\ Wang: {\it Weak multiplier Hopf algebras III. Integrals and duality}. Preprint University of Leuven (Belgium) and Southeast University of Nanjing (China) (2016). In preparation.
\snl 
{[\bf Ve]} P.\ Vecserny\'es: {\it Larson-Sweedler theorem and the role of grouplike elements in weak Hopf algebras}. J.\ Algebra 270 (2003), 471-520. See also arXiv: 0111045v3 [math.QA] for an extended version.
\snl

\end